\documentclass{mcom-preprint}
\usepackage{amscd,amssymb}
\usepackage{graphicx,color}
\copyrightinfo{2007}{D.~Arnold, G~Awanou, R.~Winther}

\numberwithin{equation}{section}
\newtheorem{theorem}{Theorem}[section]
\newtheorem{lemma}[theorem]{Lemma}

\newcommand{\ud}{\,d}
\newcommand{\Sym}{\mathbb{S}}
\newcommand{\R}{\mathbb{R}}
\newcommand\K{\mathbb{K}}
\newcommand\Mat{\mathbb{M}}
\renewcommand\P{{\mathcal P}}
\newcommand\M{{\mathcal M}}
\newcommand\N{{\mathcal N}}
\newcommand\T{{\mathcal T}}
\newcommand\tr{\operatorname{tr}}
\newcommand\skw{\operatorname{skw}}
\newcommand\sym{\operatorname{sym}}
\newcommand\x{\times}
\renewcommand{\div}{\operatorname{div}}
\newcommand{\grad}{\operatorname{grad}}
\newcommand{\eps}{\operatorname{\epsilon}}
\newcommand{\curlf}{\operatorname{curl}_f}
\newcommand{\rotf}{\operatorname{rot}_f}
\newcommand{\rotff}{\operatorname{rot}_f^*}
\newcommand{\curl}{\operatorname{curl}}
\newcommand{\gradf}{\operatorname{\grad}_f}
\newcommand{\gradff}{\operatorname{\grad}_f^* }
\newcommand{\curlff}{\operatorname{curl}_f^*}
\newcommand{\curll}{\operatorname{curl}^*}
\newcommand{\ccurl}{\operatorname{curl}\operatorname{curl}^*}
\newcommand{\rote}{\operatorname{rot}_{e^{\perp}}}
\newcommand{\XI}{\operatorname{\Xi}}
\DeclareMathOperator{\vect}{\operatorname{vec}}
\def\<{\langle}
\def\>{\rangle}
\begingroup\makeatletter\ifx\SetFigFont\undefined%
\gdef\SetFigFont#1#2#3#4#5{%
  \reset@font\fontsize{#1}{#2pt}%
  \fontfamily{#3}\fontseries{#4}\fontshape{#5}%
  \selectfont}%
\fi\endgroup%

\begin{document}
\title{Finite elements for symmetric tensors in three dimensions}

\author{Douglas Arnold}
\address{Institute for Mathematics and its Applications,
University of Minnesota,
Minneapolis, MN, 55455}
\email{arnold@ima.umn.edu}
\urladdr{http://www.ima.umn.edu/\char'176arnold}

\author{Gerard Awanou}
\address{Department of Mathematical Sciences,
Northern Illinois University,
Dekalb, IL, 60115}
\email{awanou@math.niu.edu}
\urladdr{http://www.math.niu.edu/\char'176awanou}

\author{Ragnar Winther}
\address{Centre of Mathematics for Applications and Department of
Informatics, University of Oslo, P.O. Box 1053, Blindern, 0316 Oslo, Norway}
\email{ragnar.winther@cma.uio.no}
\urladdr{http://folk.uio.no/\char'176rwinther}


\date{January 16, 2007}

\begin{abstract}
We construct finite element subspaces of the space of symmetric tensors
with square-integrable divergence on a three-dimensional domain.  These
spaces can be used to approximate the stress field in the classical
Hellinger--Reissner mixed formulation of the elasticty equations, when
standard discontinous finite element spaces are used to approximate
the displacement field.  These finite element spaces are defined
with respect to an arbitrary simplicial triangulation of the domain,
and there is one for each positive value of the polynomial degree
used for the displacements.  For each degree, these provide a stable
finite element discretization.  The construction of the spaces is closely
tied to discretizations of the elasticity complex, and can be viewed as
the three-dimensional analogue of the triangular element family for plane
elasticity previously proposed by Arnold and Winther.
\end{abstract}

\maketitle

\section{Introduction}

The classical Lagrange finite element spaces provide natural simplicial
finite element discretizations of the Sobolev space $H^1$.  Similarly, various
finite element spaces derived in the theory of mixed finite elements,
such as the Raviart--Thomas and Nedelec spaces, provide the natural finite
element discretizations of the spaces $H(\div)$ and $H(\curl)$. (These
statements are made precise and treated in a uniform framework
of the finite element exterior calculus in \cite{afw-exterior}.)
In this paper we consider the finite element discretization of the
space $H(\div,\Omega;\Sym)$ consisting of square-integrable symmetric
tensors (or, given a choice of coordinates, symmetric matrix fields)
with square-integrable divergence.  In the classical Hellinger--Reissner
mixed formulation of the elasticity equations, the stress is sought
in $H(\div,\Omega;\Sym)$  and the displacement in $L^2(\Omega;\R^n)$.
The natural discretization of the latter space is evident---piecewise
polynomial of some degree without interelement continuity
constraints---but the development of an appropriate finite element
subspace of $H(\div,\Omega;\Sym)$ to use with these is a long-standing and
challenging problem.  For plane elasticity, the known stable mixed finite
element methods have mostly involved composite elements for the stress
\cite{Arnold-Douglas-Gupta,fraejisdv,Johnson-Mercier,watwood-hartz}.
To avoid these, other authors have modified the standard mixed
variational formulation of elasticity to a formulation that uses
general, rather than symmetric, tensors for the stress, with the
symmetry imposed weakly; see \cite{Amara-Thomas,Arnold-Brezzi-Douglas,
Arnold-Falk,Stein-Rolfes,Stenberg86,Stenberg88,Stenberg-mafelap,afw-weak}.
Not until 2002 was a stable non-composite finite method for the classical
mixed formulation of plane elasticity found \cite{awi1}.  This work
can be seen as answering the question ``what are the natural finite
element discretizations of $H(\div,\Omega;\Sym)$?'' in the case of two dimensions.
In this paper we address this question for three dimensions.

Up until recently, there were no mixed finite elements for the
Hellinger--Reissner formulation in three dimensions  known to be stable.
In \cite{adams}, a partial analogue of the lowest order element in
\cite{awi1} was proposed and shown to be stable.  Here we will derive
the full analogue of the results of \cite{awi1}.  We construct a family
of finite element subspaces of $H(\div,\Omega;\Sym)$ which, when used to
discretize the stress in elasticity along with the obvious discontinous piecewise
polynomial discretization for the displacement, provide stable mixed finite
elements for the Hellinger--Reissner principle.  As in two dimensions,
these spaces lead to a finite element subcomplex of the \emph{elasticity
complex}, related to it by commuting diagrams.  The elasticity complex
reveals a close connection between the finite element discretization
of $H(\div,\Omega;\Sym)$ and finite element discretization of the space
$H(\ccurl,\Omega;\Sym)$, involving a second order operator.  A key step in
our analysis is the identifcation of the conditions which are necessary
for a piecewise smooth matrix field to belong to this space. These
conditions are derived in Section~\ref{imp_sec} below.

We recall the standard mixed formulation for the elasticity equations.
Let $\Omega$ be a contractible polyhedral domain in $\mathbb{R}^3$,
occupied by a linearly elastic body which is clamped on the boundary
$\partial \Omega$, and let $S$ and $u$ denote the stress and displacement
fields engendered by a force $f$ acting on the body.  The matrix field
$S$ and the vector field $u$ can be characterized as the unique critical
point of the Hellinger-Reissner functional
$$
\mathcal{J}(T,v) = \int_{\Omega} (\frac{1}{2} A T : T +
\div T \cdot v -
f \cdot v )  \ud x
$$
over the space $H(\div,\Omega;\Sym) \times L^2(\Omega;\mathbb{R}^3)$.
Here, $\Sym$ is the six dimensional space of symmetric matrices
and $S: T$ denotes the Frobenius product on $\Sym$.
The given
compliance tensor $A=A(x):\Sym \to \Sym$ is symmetric,
and bounded and
positive definite uniformly
with
respect to $x \in \Omega$. The divergence operator, $\div$, is applied
to a matrix field by taking the divergence of each row. Hence, this operator
maps the space $H(\div,\Omega;\Sym)$ into $L^2(\Omega;\mathbb{R}^3)$.

A mixed finite element method determines an approximate stress field
$S_h$ and
an approximate displacement field $u_h$ as the unique critical
point $(S_h,u_h)$ of the Hellinger--Reissner functional
in a finite element space
$\Sigma_h \times V_h \subset  H(\div,\Omega;\Sym) \times
L^2(\Omega;\mathbb{R}^3)$, where $h$
denotes the mesh size. Equivalently, $(S_h,u_h) \in \Sigma_h \times
V_h$
solves the saddle point system
\begin{equation}\label{system}
\int_\Omega(AS_h:T + \div T\cdot u_h +\div
S_h\cdot v)\,dx =
\int_\Omega fv\,dx, \quad
(T,v)\in \Sigma_h \times V_h.
\end{equation}
To ensure that the discrete system has a unique solution and that it
provides a good approximation of the true solution the finite dimensional
spaces
$\Sigma_h$ and $V_h$ must satisfy the stability conditions from the
theory of mixed finite element methods, see \cite{Brezzi, Brezzi-Fortin}.
As is well known, see for example \cite{awi1}, the
following two conditions are
sufficient:
\begin{itemize}
\item $\div\Sigma_h \subset V_h$.
\item There exists a linear operator $\Pi_h: H^1(\Omega;\Sym) \to
\Sigma_h$, bounded in
$\mathcal{L}(H^1;L^2)$
uniformly with respect to $h$, and such that
$\div \Pi_h \ S = \Pi^V_h \div S$
for all $S \in  H^1(\Omega;\Sym)$, where $\Pi^V_h:
L^2(\Omega;\mathbb{R}^3) \to V_h$ denotes the $L^2$-projection.
\end{itemize}
As mentioned above,
the construction of finite element spaces which fulfill these
two conditions has proved to be surprisingly hard. In this paper, we will derive
a family of finite element spaces $\Sigma_h$ and $V_h$ based
on  tetrahedral meshes, and
show that they satisfy these two
stability conditions.  There is one member of the family for
each polynomial degree $k\ge 1$.  The space $V_h$ for the displacements
is simply the space of all piecewise polynomial vector fields of
degree at most $k$.  In the lowest order case, $k=1$, the space $\Sigma_h$
contains the full space of quadratic polynomials on each element, augmented
by divergence-free polynomials of degree $3$ and $4$.  The local dimension of $\Sigma_h$ is
$162$, or $27$ per component of stress on average.  The analogous space
in two dimensions, derived in \cite{awi1}, was of local dimension $24$ ($8$
per component).

The complexity of the elements may very
well limit their practical significance. However, we believe that the determination
of the natural discretization of space $H(\div,\Omega;\Sym)$ provides important
insight, both into the obstacles to the derivation of simpler methods, and for the design
of alternative procedures, such as nonconforming methods.

This paper is organized as follows.
After giving some preliminaries remarks in Section~\ref{preliminaries}
we present the lowest
order element and establish its properties in Section~\ref{lowest-order}.
A family of higher order elements is then presented
in Section~\ref{higher-order}.
Key to the analysis of these elements is the
description of the polynomial space of symmetric matrix fields with
vanishing divergence and vanishing normal traces on the boundary of
a simplex $K$. The dimension of this space is derived in
Section~\ref{basis}
based on preliminary results derived in Sections~\ref{imp_sec}
and \ref{single-tet}. Furthermore, an explicit basis for this space,
necessary
for the computational procedure, is also given in the lowest order cases.
Finally, in Section~\ref{elasticity-complex} we summarize the
results of our construction by
presenting a discrete analogue of the elasticity complex.

\section{Notation and preliminaries} \label{preliminaries}
We begin with some basic notation.
If $K \subset \R^3$ is a tetrahedron, then $\Delta_2(K)$ denotes the
set of the four $2$--dimensional faces of $K$, $\Delta_1(K)$
the set of the six $1$--dimensional edges, and $\Delta_0(K)$
the set of the four vertices. Furthermore, $\Delta(K)$ is the set of
all subsimplexes of $K$ (of dimension $0,1,2$ or $3$).

We let $\Mat$ be
the space of $3\times 3$ real matrices, and
$\Sym$ and $\K$ the subspaces of symmetric and skew symmetric matrices,
respectively.  The operators $\sym:\Mat\to\Sym$ and $\skw:\Mat\to\K$ denote
the symmetric and skew symmetric parts, respectively.
Note that an element of the space $\K$ can be identified with its
\emph{axial vector} in $\R^3$ given by the map
$\vect : \K \to \R^3$:
\begin{equation*}
\vect \begin{pmatrix} 0 & -v_3 & v_2 \\
v_3 & 0 & -v_1 \\ -v_2 & v_1 & 0 \end{pmatrix}
= \begin{pmatrix} v_1 \\ v_2 \\ v_3 \end{pmatrix},
\end{equation*}
i.e., $\vect^{-1}(v)w= v\x w$ for any vectors $v$ and $w$.

For any vector space $X$, we let $L^2(\Omega;X)$ be the
space of square-integrable vector fields on $\Omega$ with values in
$X$. For our
purposes, $X$ will usually
either be $\R$,  $\mathbb{R}^3$, or $\Mat$, or some subspace
of one of these. In the case $X=\mathbb{R}$, we will simply write $L^2(\Omega)$.
The corresponding Sobolev space of order $k$,
i.e., the subspace of $L^2(\Omega;X)$  consisting of functions with
all partial derivatives of order less than or equal to $k$
in $L^2(\Omega;X)$, is denoted $H^k(\Omega;X)$, and its norm by
$\| \,\cdot \,\|_k$.
The space $H(\div,\Omega;\Sym)$ is defined by
\[
H(\div,\Omega;\Sym) = \{ T \in L^2(\Omega;\Sym)\, | \, \div T \in
L^2(\Omega;\R^3) \},
\]
where the divergence of a matrix field is ithe vector field obtained by
applying the
divergence operator row--wise.
For a vector field $v: \Omega \to \mathbb{R}^3$, $\grad v$ is the matrix field
with rows the gradient of each component, and the symmetric gradient,
$\eps(v)$, is given by
$\eps(v) = \sym \grad v$. Furthermore,
\[
\curl v = - 2\vect \skw \grad v =
\left(\begin{array}{c} \partial_3 v_2 - \partial_2 v_3 \\
-\partial_3 v_1 + \partial_1 v_3\\
\partial_2 v_1 - \partial_1 v_2
\end{array} \right).
\]
If we consider a linear coordinate transformation of the form $x = B
\hat x + b$, with corresponding vector fields $v$ and $\hat v$ related
by
$v(x) = (B')^{-1}\hat v (\hat x)$, then we have
\begin{equation*}
\skw \grad_x v = (B')^{-1}(\skw \grad_{\hat x} \hat v )B^{-1}.
\end{equation*}
Here $B'$ denotes the transpose of $B$.
In particular, if $B$ is orthogonal, i.e., $B'B = I$, then $v = B\hat
v$
and
\begin{equation}\label{transformation_o}
\skw \grad_x v = B (\skw \grad_{\hat x} \hat v )B'.
\end{equation}

As for the divergence and the gradient operator, the operator
$\curl$ acts on a matrix field by applying the ordinary curl operator
to each row of the matrix.
The operator $\curl^*$ is the corresponding
operator obtained by taking the curl of each column. Alternatively, we
have $\curl^* T = (\curl T')'$ for any matrix field T.
The second order operator $\ccurl$ maps symmetric matrix fields
into symmetric matrix fields. Let $\XI :\Mat \to \Mat$ be the algebraic operator
$\XI T = T' - \tr(T)I$,
where $I$ is the identity matrix. Then $\XI$ is invertible with
$\XI^{-1}T = T' - \tr(T)I/2$. The following identities are useful:
\begin{align}\label{skewcurl}
\vect \skw \curl T &= -\frac{1}{2} \div \XI T, \quad T \in
C^\infty(\Omega, \Mat),\\ \label{curlskew}
\curl T &= \XI \grad \vect T, \quad T \in
C^\infty(\Omega, \K),\\ \label{tracecurl}
\tr \curl T &= - 2 \div \vect \skw T, \quad T \in
C^\infty(\Omega, \Mat).
\end{align}
These formulas can be verified directly, but they are also consequences of
the discussions given in \cite{afw-weak, afw-exterior}, cf.~Section~4
of \cite{afw-weak} or Section~11 of \cite{afw-exterior}.

For $K \subset \R^3$ we let
$\mathcal{P}_{k}(K;X)$ be the space of polynomials of degree $k$,
defined on $K$ and with values in $X$. We write $\mathcal{P}_{k}$
or $\mathcal{P}_k(K)$
for $\mathcal{P}_{k}(K;\mathbb{R})$. The de Rham complex has
a polynomial analogue of the form
\begin{equation}\label{deRham}
\R \hookrightarrow \P_{k+3} \xrightarrow{\grad} \P_{k+2}(K;\R^3) \xrightarrow{ \curl} \P_{k+1}(K;\R^3)
\xrightarrow{\div } \P_{k} \to 0.
\end{equation}
In fact, this complex is an exact sequence \cite{afw-exterior}.

In recent years differential complexes have come to play a
significant role in the design of mixed finite
element methods  \cite{a,awi1,afw-weak,afw-exterior}. For the
equations of elasticity, the relevant differential
complex is the \emph{elasticity complex}.  In three space
dimensions,
the elasticity complex takes the form
\begin{equation*}
\mathbb{T} \hookrightarrow C^{\infty}(\Omega;\mathbb{R}^3) \xrightarrow{\eps} C^{\infty}(\Omega;\Sym) \xrightarrow{\ccurl} C^{\infty}(\Omega;\Sym)
\xrightarrow{\operatorname{div} } C^{\infty}(\Omega;\mathbb{R}^3)  \to 0,
\end{equation*}
where $\mathbb{T}$ is the six-dimensional space of infinitesimal rigid
motions, i.e., the space of linear polynomial functions of the form
$x\mapsto a+b \times x$ for some $a,b \in \mathbb{R}^3$. It is straightforward to verify that the
elasticity complex is a complex, i.e., the composition of two
successive
operators is zero. In fact, if the domain $\Omega$ is contractible, then the
elasticity complex is an exact sequence; see \cite{afw-weak,afw-exterior}.

An analogous complex with less smoothness is
\begin{equation*}
\mathbb{T} \hookrightarrow H^{1}(\Omega;\mathbb{R}^3)
\xrightarrow{\eps} H( \ccurl,\Omega;\Sym)
\xrightarrow{\ccurl} H(\div,\Omega;\Sym)
\xrightarrow{\div} L^2(\Omega;\mathbb{R}^3)  \to 0,
\end{equation*}
where $ H( \ccurl,\Omega;\Sym)=\{ \, S \in L^2(\Omega;\Sym) \, | \, \ccurl S \in L^2(\Omega;\Sym) \, \}$.

There is also a polynomial analogue of the elasticity complex. Let
$K \subset \R^3$ be tetrahedron and $k \geq 0$. The polynomial
elasticity complex  is given by
\begin{equation}\label{cd1}
\mathbb{T} \hookrightarrow
\mathcal{P}_{k+4}(K;\mathbb{R}^3) \xrightarrow{\eps} \mathcal{P}_{k+3}(K;\Sym) \xrightarrow{\ccurl}
\mathcal{P}_{k+1}(K;\Sym) \xrightarrow{\div }
\mathcal{P}_{k}(K;\mathbb{R}^3)\to 0.
\end{equation}
This complex is an exact sequence.
To prove the exactness, we first show that if $S$ is a matrix field in $\mathcal{P}_{k+3}(\Omega;\Sym)$, with
$\ccurl S =0$, then $S=\eps(u)$ for $u=(u_1,u_2,u_3)^T \in \mathcal{P}_{k+4}(\Omega;\mathbb{R}^3)$. Clearly
$S=\eps(u)$ for $u \in  C^{\infty}(K;\mathbb{R}^3)$. It is enough to show that all second derivatives of $u_k$, $k=1$, $2$, $3$, are
in $\mathcal{P}_{k+2}(\Omega;\mathbb{R})$.
This follows from the identity
$\partial_{ij} u_k = \partial_i \eps(u)_{jk} + \partial_j \eps(u)_{ik} - \partial_k \eps(u)_{ij}$.
Next, we show that if $S \in \mathcal{P}_{k+1}(\Omega;\Sym)$, and $\div S =0$, then $S=\ccurl T$ for some $T \in
\mathcal{P}_{k+3}(\Omega;\Sym)$.
First we observe that since $\div S =0$ it follows from the fact that
\eqref{deRham} is exact that $S = \curl U$ for some $U \in
\P_{k+2}(K;\Mat)$.
Furthermore, since $S$ is symmetric it follows from \eqref{skewcurl}
that $\div \XI U = 0$, and as a consequence, using \eqref{deRham}
once more, we obtain that $\XI U = \curl T$ for some $T \in
\P_{k+3}(K;\Mat)$,
or
$S = \curl U = \curl \XI^{-1}\curl T$.
However, by \eqref{curlskew} we have
$\curl \XI^{-1}\curl \skw T = \curl \grad \vect \skw T = 0$.
Hence, we can take $T \in \P_{k+3}(K;\Sym)$. Finally, we observe that
if $T$ is symmetric, then \eqref{tracecurl} implies that $\tr \curl T =
0$,
and therefore
$S = \curl \XI^{-1}\curl T = \ccurl T$.
To establish the surjectivity of the last map, one can use the fact
that
$\dim \P_k = (k+1)(k+2)(k+3)/6$ to verify that the
alternating sum of the dimensions
of the spaces in the sequence is zero.
The same arguments show that \eqref{cd1} is exact for $k=-1$, $-2$, or
$-3$, if $\mathcal{P}_j$ is interpreted as the zero space for $j<0$.

Let $\{\mathcal{T}_h\}$ denote a family of triangulations of $\Omega$
by tetrahedra with diameter bounded by $h$. We assume that the intersection of
any two tetrahedra in $\T_h$
is either empty or a common subsimplex of each. The
family $\{\mathcal{T}_h\}$ is also assumed to be
shape regular in the sense that the
ratio of the radii of the circumscribed and inscribed spheres
of all the tetrahedra can be bounded by a fixed constant.
Furthermore, we will use the notation $\Delta_j(\T_h)$, for $j=0,1,2,$
to denote the set of vertices, edges, and faces, respectively,
associated with the mesh $\T_h$.
In Section~\ref{higher-order} we will define a family
of finite element spaces $\Sigma_h \subset H(\div,\Omega;\Sym$)
and $V_h \subset L^2(\Omega; \mathbb{R})$ for the elasticity problem
consisting of piecewise polynomial spaces with repect to $\mathcal{T}_h$
of arbitrarily
high polynomial order.
However, we will first consider the lowest order case of this family
in Section~\ref{lowest-order} below. All our spaces will have the
property that
$\div \Sigma_h \subset V_h$.
Furthermore, we will identify a corresponding projection operator
$\Pi_h:H^1(\Omega;\Sym)\to \Sigma_h$ satisfying the commutativity
relation
\begin{equation}\label{com1}
\div \Pi_h T = \Pi^V_h \div  T, \qquad T \in H^1(\Omega;\Sym),
\end{equation}
and the bound
\begin{equation}\label{nind}
\|\Pi_h T\|_0 \leq C \|T\|_1, \qquad T \in H^1(\Omega;\Sym),
\end{equation}
with constant $C$ independent of $h$. Here $\Pi^V_h : L^2(\Omega;\R^3) \to
V_h$
is the $L^2$ projection.
It is a consequence of the general error bounds derived in
\cite{Falk-Osborn}, cf.~also \cite{awi1}, that the properties above
imply that $(\Sigma_h,V_h)$ is a stable pair of elements for the
discretization \eqref{system}, and that the error bounds
\begin{gather}\label{S-bound}
\|S - S_h\|_0 \le \|(I - \Pi_h)S \|_0\\ \label{u-bound}
\| u - u_h \|_0 \le \| (I - \Pi^V_h)u \|_0 + c\|(I - \Pi_h)S \|_0
\end{gather}
holds, with a constant $c$ independent of $h$.
Here $(S,u)$ is the unique critical point of the Hellinger-Reissner
functional
over $H(\div,\Omega;\Sym) \times L^2(\Omega;\mathbb{R}^2)$ and
$(S_h,u_h) \in \Sigma_h \x V_h$ the corresponding finite element solution.
In addition, $\div S_h = \Pi^V_h \div S$.

\section{The lowest order element} \label{lowest-order}
We first describe the restriction of the lowest order
spaces $\Sigma_h$ and $V_h$
to a single tetrahedron $K \in \mathcal T_h$.
Define
\begin{equation*}
\Sigma_K = \{\, T \in \mathcal{P}_4(K;\mathbb{S}) \, | \,
 \div T \in \mathcal{P}_1(K;\mathbb{R}^3) \, \}, \quad
V_K = \mathcal{P}_1(K;\mathbb{R}^3).
\end{equation*}
The space $V_K$ has dimension 12 and a complete set of degrees of
freedom
is given by the zero and first order moments with respect to $K$.
The space $\Sigma_K$ has dimension at least 162 since the dimension of
$\mathcal{P}_4(K;\mathbb{S})$ is 210 and the condition $\div T \in \mathcal{P}_1(K;\mathbb{R}^3)$
represents 48 linear constraints.  We will show that $\dim \Sigma_K =
162$
by
exhibiting $162$ degrees of freedom which determine the elements uniquely.
Define
$$
\M(K)=\{ \, T \in \mathcal{P}_{4}(K;\Sym) \, | \,
 \div T=0,\ T  n=0 \text{\quad on $\partial K$} \, \}.
$$
It will follow from Theorem~\ref{dimension-formula} below
that the
dimension of $\M(K)$ is 6. Furthermore, in
Section~\ref{basis} we will also
give an explicit basis for this space.
Using this basis, we can state the $162$ degrees of freedom for the
space
$\Sigma_K$.
\begin{lemma}\label{lem1}
A matrix field $T \in \Sigma_K$ is uniquely determined by the following
degrees of freedom
\begin{enumerate}
\item the values of $T$ at the vertices of $K$, $4 \times 6=24$
degrees of freedom,
\item for each edge $e \in \Delta_1(K)$ with unit tangent vector $s$
and linearly independent normal vectors $n_-$ and $n_+$,
the  constant, linear
and quadratic moments over $e$ of $s'T n_-$,  $s'T n_+$,
$n'_-T n_-$,  $n'_+T n_+$,  $n'_-T n_+$, $6
\times 3 \times 5=90$ degrees of freedom,
\item for each face $f \in \Delta_2(K)$, with normal $n$, the constant
  and
linear moments over $f$ of $T n$, $4 \times 3 \times 3=36$ degrees of freedom,
\item the average of $T$ over $K$, 6 degrees of freedom,
\item the value of the moments $\int_K T:U \ \ud x$, $U \in
  \M(K)$,
$6$ degrees of freedom.
\end{enumerate}
\end{lemma}
\begin{proof}
We assume that all degrees of freedom vanish and show that $T=0$.
Since $T=0$ at the vertices, the second set of degrees of freedom
imply that $T n=0$ on each edge for both faces meeting the edge.
By the third set of degrees of freedom we obtain that
$T n=0$ on each face of $K$.
For
$v=\div T \in \P_1(K,\R^3)$ we have
$$
\int_K v^2 \ud x = - \int_K T:\eps(v) \ud x + \int_{\partial K} T n
\cdot v \ud x_f = - \int_K T:\eps v \ud x = 0
$$
by the fourth set of degrees of freedom.
Here and below, $d x_f$ denotes the surface measure on $\partial K$.
We conclude that $\div T=0$,
and,
by
the last set of degrees of freedom, that $T =0$.
\end{proof}

We now describe the finite element spaces on the triangulation $\mathcal{T}_h$. We
denote by $V_h$
the space
of vector fields which belong to $V_K$ for each $K \in
\mathcal{T}_h$ and by
$\Sigma_h$ the space of matrix fields which belong piecewise to
$\Sigma_K$, and with the continuity conditions induced by the degrees
of freedom. In particular, for $T\in\Sigma_h$, the normal components
$Tn$ are continuous across all faces
$f \in \Delta_2(\mathcal{T}_h)$ and,
hence, $\Sigma_h \subset H(\div,\Omega;\Sym)$.  In addition, 
if $T\in\Sigma_h$, $e \in \Delta_1(\T_h)$, and $s$ and $n$ are vectors which are tangential
and normal to $e$, respectively, then $s'Tn$ is continuous on $e$.

It remains to define an interpolation operator
$\Pi_h : H^1(\Omega;\Sym) \to \Sigma_h$ which satisfies
\eqref{com1} and \eqref{nind}.
Because of the vertex and edge
degrees of freedom, the canonical interpolation operator for
$\Sigma_h$, $\Pi_h^{\Sigma}$, defined directly from
the degrees of freedom, is not
bounded on
$ H^1(\Omega; \Sym)$. In order to overcome this difficulty we
introduce the operator $\Pi_h^0 : H^1(\Omega; \Sym) \to \Sigma_h$
defined from the degrees of freedom above, but where
the vertex and edge degrees of freedom are set equal to zero,
i.e., we have
\begin{align}
\Pi^0_h T(x) & = 0,  \quad x \in \Delta_0(\T_h),
\label{i1}
\\
\int_e \Pi^0_h T n \cdot v \ud s & = 0,
\quad e \in \Delta_1(\T_h),\, v \in \mathcal{P}_2(e;\R^3),\,
n\in e^\perp, \label{i2}
\\
\int_f (T - \Pi_h^0 T) n \cdot v \ud x_f & = 0, \quad f \in \Delta_2(\T_h),
\, v \in \mathcal{P}_1(f;\R^3),  \label{i3} \\
\int_K (T - \Pi_h^0 T) \ud x & = 0,  \quad K \in T_h,  \label{i4} \\
\int_K (T - \Pi_h^0 T):U \ud x & = 0,  \quad K \in T_h, \,
U \in \M(K). \label{i5}
\end{align}
The commutativity property \eqref{com1} for $\Pi_h^0$ follows from
\eqref{i3} and \eqref{i4} since
$$
\int_K (\div \Pi^0_h T - \div T) \cdot v \ud x =
- \int_K (\Pi^0_h T - T): \eps(v) \ud x +
\int_{\partial K} (\Pi^0_h T - T)n \cdot v \ud x_f = 0.
$$
The uniform boundedness \eqref{nind} can be seen from a standard
scaling argument using the matrix Piola transform.
Let $\hat K$ be a fixed reference tetrahedron and $F=F_K :\hat{K} \to
K$ be
an affine isomorphism of the form $F \hat{x} = B \hat{x} +b$. Given
a matrix field $\hat T: \hat{K} \to \mathbb{S}$, define $T: K \to
\mathbb{S}$
by the matrix Piola transform
$T(x) = B \hat{T}(\hat{x}) B^T$,
with $x=F\hat{x}$. Using $\div T(x) = B \div \hat{T}(\hat{x})$,
it is easy to verify that
$T \in \Sigma_K$ if and only if $\hat T \in \Sigma_{\hat{K}}$.
Furthermore, as in \cite{arw,awi1} a scaling
argument can be used
to verify the uniform boundedness condition
\eqref{nind} for the operator $\Pi_h^0$.
We can therefore conclude that the operator $\Pi_h^0$ satisfies the
two conditions \eqref{com1} and \eqref{nind}.
However, the operator $\Pi_h^0$ lacks good approximation properties.
Therefore, in order to obtain error estimates from the general bounds
\eqref{S-bound} and \eqref{u-bound} a more accurate interpolation
operator is needed.

Consider the modified interpolation operator
$\Pi_h : H^1(\Omega;\Sym) \to \Sigma_h$ of the form
\begin{equation}\label{opdef}
\Pi_h = \Pi^0_h (I - R_h) + R_h,
\end{equation}
where $R_h : L^2(\Omega;\Sym) \to \Sigma_h$ is the Cl\'ement operator
onto the continuous piecewise quadratic subspace of $\Sigma_h$
\cite{c}.
This operator satifies the bounds
\begin{equation*}
\|R_h T - T \|_j \leq c h^{m-j} \|T\|_m, \quad 0\leq j \leq 1,
\quad j \leq m
\leq 3.
\end{equation*}
As a consequence of this bound, and the boundedness \eqref{nind} of
$\Pi_h^0$,
we obtain the estimate
\begin{equation}\label{Pi-estimate}
\|\Pi_h T - T \|_0 \leq c h^m \|T\|_m, \quad 1 \leq m \leq 3
\end{equation}
for the interpolation error. Furthermore, since $\Pi_h^0$ satisfies
\eqref{com1} and $\Pi^V_h \div R_h = \div R_h$, we conclude that
$\Pi_h$ satisfies \eqref{com1}.

We also recall that the projection operator $\Pi^V_h : L^2(\Omega;\R^3)
\to V_h$
satisfies the error estimate
\begin{equation}\label{P-estimate}
\|\Pi^V_h v - v \|_0 \leq c h^m \|v\|_m,  \quad 0 \leq m \leq 2.
\end{equation}
The estimates \eqref{Pi-estimate} and \eqref{P-estimate}, combined
with the basic error bounds
\eqref{S-bound} and \eqref{u-bound}, and the fact that $\div S_h =
\Pi^V_h\div S$, imply the following error estimates for the finite
element method
generated by $\Sigma_h \x V_h$.

\begin{theorem}\label{convergence-lowest order}
Let ($S,u$) denote the unique critical point of the Hellinger-Reissner
functional over $H(\div,\Omega;\Sym) \times L^2(\Omega;\mathbb{R}^2)$
and let ($S_h,u_h$) be the unique critical point over $\Sigma_h \times
V_h$. Then
\begin{eqnarray*}
\|S -S_h \|_0 & \leq & c h^m \|S\|_m, \qquad 1\leq m \leq 3,\\
\| \div S -\div S_h \|_0 &
\leq & c h^m \|\div S \|_m, \qquad 0\leq m \leq 2,\\
\|u-u_h \|_0 & \leq & c h^m \|u \|_{m+1}, \qquad 1\leq m \leq 2.
\end{eqnarray*}
\end{theorem}

{\it Remark.} In general it seems not possible to lower the polynomial degree
of the stress space $\Sigma_h$ introduced above. However, as in
\cite{awi1},
a minor simplification is possible. On each tetrahedron $K \in
\mathcal T_h$ we take the restricted displacement space
$V_K$ to be the rigid motions $\mathbb{T} \subset \P_1(K;\R^3)$
and the corresponding stress space to be
$$
\tilde \Sigma_K = \{\, T \in \mathcal{P}_{4}(K;\Sym) \, | \, \div T \in \mathbb{T} \, \}.
$$
Clearly dim $\tilde \Sigma_K \geq 210-(60-6)=156$. In fact, $\dim
\tilde \Sigma_K =
156$
and a complete set of degrees of freedom is obtained by removing
the six average values of $T$ represented by (4) in Lemma~\ref{lem1}.
The proof of the fact that these degrees of freedom are unisolvent
for $\tilde \Sigma_K$ follows by a simple modification of the proof of
Lemma~\ref{lem1} above. Just observe that if $v = \div T \in
\mathbb{T}$
then $\eps(v) = 0$. However, the simplified element is less
accurate,
since the stress space lacks some quadratics,
and the displacement space some linears.
Instead of the error estimates given in
Theorem~\ref{convergence-lowest order}
we obtain at most $O(h^2)$ convergence for $||S -S_h ||_0$,
and at most first order convergence for $|| \div (S - S_h) ||_0$
and $||u-u_h ||_0$.

\section{A Family of Higher Order Elements} \label{higher-order}
In this section we describe a family of stable element pairs, one for
each degree $k\ge 1$.  The lowest order case $k=1$ is the one treated above.
We first describe the elements  on a single tetrahedron. Define
\begin{equation*}
\Sigma_K = \{\, T \in \mathcal{P}_{k+3}(K;\Sym) \, | \, \div T \in \mathcal{P}_k(K;\mathbb{R}^3) \, \}, \quad
V_K = \mathcal{P}_k(K;\mathbb{R}^3).
\end{equation*}
Then
$$
\dim V_K = 3 \binom{k+3}{3} =   \frac{(k+3)(k+2)(k+1)}{2} ,
$$
\begin{align*}
\dim \Sigma_K \geq d_k  :&= \dim \mathcal{P}_{k+3}(K;\mathbb{S})-[\dim
\mathcal{P}_{k+2}(T;\mathbb{R}^3)- \dim \mathcal{P}_{k}(T;\mathbb{R}^3)]\\
&=6 \binom{k+6}{3} - 3 \binom{k+5}{3} + 3 \binom{k+3}{3} = k^3+12k^2+56k+93.
\end{align*}
Notice that the space $\eps[\mathcal{P}_k(K,\mathbb{R}^3)]$ has dimension
$(k+3)(k+2)(k+1)/2-6$. Analoguous to the lowest order case,
we define the space
\begin{equation*}
\M_k(K)=\{\, T \in \mathcal{P}_{k}(K;\Sym)\, |
\, \div  T=0,\quad T  n=0 \text{\quad on $\partial K$} \, \}.
\end{equation*}
We will prove in Section~\ref{basis}, Theorem~\ref{dimension-formula},
that $\dim\M_k(K)$ is $(k+2)(k-2)(k-3)/2$
for $k \ge 4$.

The degrees of freedom for $V_K$ are the moments of degree less than or
equal to $k$ with respect to $K$.
A unisolvent set of degrees of freedom for $\Sigma_K$ are given by
\begin{enumerate}
\item the values of $T$ at the vertices of $K$, $4 \times 6=24$
  degrees of freedom,
\item for each edge $e \in \Delta_1(K)$ with unit tangent vector $s$
and linearly independent normal vectors $n_-$ and $n_+$,
the moments
of degree at most $k+1$ over $e$ of $s'T n_-$,  $s'T n_+$,
$n'_-T n_-$,  $n'_+T n_+$,  $n'_-T n_+$,
$6 \times (k+2) \times 5=30k+60$ degrees of freedom,
\item for each $f \in \Delta_2(K) $ with normal vector $n$, the
moments of degree at most $k$ over $f$ for $Tn$,
 $3 \times 4 \times (k+2)(k+1)/2=6 k^2+18k+12$ degrees of freedom,
\item   $\int_K T: U \ud x$, $U \in \eps(V_K)$,
$(k+3)(k+2)(k+1)/2-6$ degrees of freedom,
\item $\int_K T: U \ud x$, $U \in \M_{k+3}(K)$,
$(k+5)(k+1)k/2$ degrees of freedom.
\end{enumerate}
The proof that this set of functionals is unisolvent for the space
$\Sigma_K$
is almost identical to the lowest order case, and it is easily
checked that their numbers sum up to $d_k$. Hence, we have shown that $\dim
\Sigma_k = d_k$. Furthermore, in Section~\ref{basis}  we will give an explicit
basis for the space $\M_k(K)$ when $k=4$ and
$k=5$.

The finite element space $V_h \subset L^2(\Omega;\R^3)$
consists of all vector fields which belong to $\P_k(K;\R^3)$
for each $K \in \T_h$, while the corresponding stress space
$\Sigma_h$ is the space of matrix fields which belong piecewise to
$\Sigma_K$, and with the continuity conditions induced by the degrees
of freedom. In particular, this implies that the normal components
$Tn$, for $T \in \Sigma_h$,  are continuous over all faces in
$\Delta_2(\T_h)$.
Hence, as in the lowest order case we have that
$\Sigma_h \subset H(\div,\Omega;\Sym)$.

The $L^2$ projection $\Pi^V_h$ onto $V_h$ satisfies the estimate
\begin{equation}\label{P-estimate2}
\|\Pi^V_h v - v \|_0 \leq c h^m \|v\|_m,  \quad 0 \leq m \leq k+1.
\end{equation}
We also
introduce the Cl\'ement interpolant $R_h : L^2(\Omega;\Sym) \to
\Sigma_h$ defined as the $\Sym$--valued
version of the standard scalar Cl\'ement interpolant into continuous
piecewise polynomials of order $k+1$. Hence, the operator $R_h$ satisfies
$$
\|R_h T -T\|_j \leq c h^{m-j} \|T\|_m, \quad 0 \leq j\leq 1, \ j \leq m \leq
k+2.
$$
Furthermore, we define the modified canonical interpolation
operator $\Pi_h$ by \eqref{opdef},
where the operator $\Pi_h^0$ is defined in complete analogy with the
lowest order case, by setting the degrees freedom associated with
$\Delta_0(\T_h)$ and $\Delta_1(\T_h)$ equal to zero.
Then the operator $\Pi_h$ satisfies
\eqref{com1} and \eqref{nind}, and the error bound
\begin{equation}\label{Pi-estimate2}
\|\Pi_h T - T \|_0 \leq c h^m \|T\|_m, \quad 1 \leq m \leq k+2.
\end{equation}
As above,
the interpolation estimates \eqref{P-estimate2} and
\eqref{Pi-estimate2},
and the error bounds
\eqref{S-bound} and \eqref{u-bound}, leads to the following
error estimates.

\begin{theorem}
Let ($S,u$) denote the unique critical point of the Hellinger-Reissner
funtional over $H(\div,\Omega;\Sym) \times L^2(\Omega;\mathbb{R}^2)$
and let ($S_h,u_h$) be the unique critical point over $\Sigma_h \times
V_h$. Then
\begin{gather*}
\|S -S_h \|_0  \leq  c h^m \|S\|_m, \qquad 1\leq m \leq k+2,\\
\| \div S -\div S_h \|_0
\leq  c h^m \|\div S \|_m, \qquad 0\leq m \leq k+1,\\
\|u-u_h \|_0  \leq  c h^m \|u \|_{m+1}, \qquad 1\leq m \leq k+1.
\end{gather*}
\end{theorem}

\section{Some properties of vector fields and matrix fields}
\label{imp_sec}
In order to perform a more detailed analysis of the spaces $\mathcal
M_k(K)$ and the finite element spaces $\Sigma_h$ we will need
some basic properties of vector fields and matrix fields. These
properties will be reviewed in the present section.

\subsection{Identities for vector fields} \label{prelim}

Let $n$ denote a fixed unit vector in $\mathbb{R}^3$,
$P_n =nn'$ the orthogonal projection onto $\R n$,
$f=n^\perp$ the plane orthogonal to $n$, and
$Q_n=I-P_n$ is the orthogonal projection onto $f$. Furthermore, set
$$
C_n = \left(\begin{array}{ccc} 0 & n_3 & -n_2 \\
-n_3 & 0 & n_1 \\
n_2 &-n_1 &0
\end{array} \right) = - \vect^{-1} n,
$$
so that $C_n v = v \times n$.
The following identities are easily checked:
\begin{equation*}
C_n'  = - C_n,\quad C_n^2  =-Q_n,  \quad C_n P_n =0,  \quad C_n Q_n = Q_n C_n =C_n.
\end{equation*}
For any vector field $v=(v_1,v_2,v_3)'$ in $\mathbf{R}^3$,
we obviously have
$v = P_n v + Q_n v$ and
\[
\curl v = \curl P_n v + \curl Q_n v
= P_n \curl P_n v + P_n \curl Q_n v + Q_n \curl P_n v + Q_n \curl Q_n v .
\]
It is elementary to verify that
\begin{equation*}
P_n \curl P_n v=0,\quad
Q_n \curl Q_n v = C_n \frac{\partial v}{\partial n},
\end{equation*}
if $n = e_3 = (0,0,1)'$.
In view of the transformation formula \eqref{transformation_o},
these identities hold for an arbitrary unit vector $n$.
Furthermore, we define
\begin{equation}\label{rot-div}
\rotf v = P_n \curl v = P_n \curl Q_n v = - (\div C_nv) n,
\quad \curlf v = Q_n \curl P_n v.
\end{equation}
With this notation we obtain the decomposition
\begin{equation}\label{curlrel}
\curl v = \rotf v + \curlf v + C_n \frac{\partial v}{\partial n}.
\end{equation}
We also define the tangential gradient $\gradf \phi =Q_n\grad \phi$,
for a scalar field $\phi$.
For a vector field $v$, $\gradf v = (\grad v)Q_n$ is the matrix field
with rows equal to the tangential gradients
of the components of $v$, and we let
\[
\eps_f(v) = \frac{1}{2} \{\gradf (Q_n v) + [\gradf(Q_n v)]'\} =Q_n \eps(v) Q_n
\]
be the tangential part of the symmetric gradient.
Note that
the definitions of $\gradf v$, $\eps_f (v)$, $\curlf v$, and $\rotf
v$
do not depend on the choice of the unit normal $n$ to $f$. The identity
\begin{equation}\label{rel01}
\curlf v = -C_n \ \grad (n'v) = -C_n \ \gradf (n'v)
\end{equation}
can be easily verified in the special case $n = e_3$, and
holds in general.

\subsection{Identities for matrix fields}
We extend the operators curl, $\curlf$ and $\rotf$ to act on and yield $3 \times 3$ matrix fields
by applying the vector operations row-wise.
More precisely, $\rotf S = (\curl S) P_n = (\curl S Q_n) P_n$
and $\curlf S = (\curl SP_n) Q_n$.
We notice that, for any constant matrix $A$,
$\curl A  S = A \curl S$.
We also recall that
$\curll S = (\curl S')'$ is the corresponding operator obtained by
applying
the curl operation to each
column. The corresponding column operators $\rotff$ and $\curlff$
are defined similarly, i.e., $\rotff S = P_n \curll S =P_n \curll Q_n
S$
and  $\curlff S = Q_n \curll P_n S.$
For a given row vector $v$,
$\gradff v = (\gradf v')' = Q_n (\grad v')'$, which is the matrix whose
columns are the tangential gradients of the components of $v$.

We now extend the decomposition \eqref{curlrel} to $\curl$ and $\curll.$
It is easy to see that $C_n S$ results in $C_n$
applied to the columns of $S$, while $-S C_n$ is $C_n$ applied
row-wise.
It follows that
\begin{equation} \label{rel05}
\curl S = \curlf S + \rotf S - \frac{\partial S}{\partial n} C_n,
\end{equation}
and
\begin{equation}\label{rel06}
\curll S = \curlff S + \rotff S + C_n \frac{\partial S}{\partial n},
\end{equation}
where $\partial S / \partial n$ is obtained by taking the directional
derivative of each component. Furthermore, the identities
\begin{equation}\label{rel7}
\curlf S = (\gradf S n) C_n, \quad
\curlff S = -C_n \ \gradff n'S
\end{equation}
are just matrix analogues of \eqref{rel01}. Note also that from the
definitions of the operators $\rotf$ and $\rotff$ we have
\begin{align} \label{rel12}
P_n (\curl \curll S) P_n  &= P_n (\rotf \curll S)
=  \rotf ( P_n \curll S)
\\
&= \rotf \rotff S = \rotf \rotff (Q_n S Q_n).
\end{align}
 We will also need exact sequences relating spaces of functions
defined on a two dimensional space. Let $f=n^\perp$. In analogy with \eqref{cd1} the following
two dimensional complexes are exact:
\begin{equation} \label{pol-elas-2d1}
\mathbb{T}_f \hookrightarrow \mathcal{P}_{k+3}(f;Q_n\R^3)
\xrightarrow{\eps_f} \mathcal{P}_{k+2}(f;Q_n\Sym Q_n)
\xrightarrow{\rotf \rotff}
\mathcal{P}_{k}(f;\R P_n)
\to 0,
\end{equation}
\begin{equation} \label{pol-elas-2d2}
\P_1(f;\R) \hookrightarrow \mathcal{P}_{k+3}(f;\mathbb{R}) \xrightarrow{\gradf \gradff} \mathcal{P}_{k+1}(f;Q_n\Sym Q_n)
\xrightarrow{\rotf}
\mathcal{P}_{k}(f;Q_n\Sym P_n)
\to 0.
\end{equation}
Here, $\mathbb{T}_f$ is
the 3-dimensional space of vector fields on $f$ of the form
$v(x)=Q_n w(x)$ for some $w\in\mathbb{T}$.

\subsection{Integration by parts}
Above we discussed the operators $P_n$, $Q_n$, and $C_n$ with respect
to
a fixed linear space $f$ with a unit normal vector $n$.
If $\Omega$ is a bounded subset of $\mathbb{R}^3$, with a piecewise
smooth boundary $\partial \Omega$,
we define these operators on $\partial \Omega$ with respect to the
tangent space $f = f(x)$ and the outer unit normal vector $n=n(x)$
for $x \in \partial \Omega$. Furthermore, the
tangential differential operators $\gradf$, $\curlf$, $\rotf$,
$\eps_f$
and so on are then defined pointwise on $\partial \Omega$.

If $v$ is a smooth vector field on $\bar{\Omega}$, then
\begin{equation*}
\int_{\Omega} v \cdot \curl w \ \ud x = \int_{\Omega} \curl v \cdot w \ \ud x
- \int_{\partial \Omega} C_n v \cdot w \ \ud x_f,
\end{equation*}
where $n$ denotes the outward unit normal to $\Omega$, and
$dx_f$ is the surface measure on $\partial \Omega$.
Moreover, if $v$ and $w$ are smooth
vector fields on $\partial \Omega$, and $f$ is a piecewise smooth
submanifold, then we have
\begin{equation*}
\int_{f} \rotf v \cdot w \ \ud x_f = \int_{f} v \cdot \curlf w \ \ud x_f
- \int_{\partial f} (s'v) (n'w) \ \ud s,
\end{equation*}
where $d s$ denotes the arc length measure on $\partial f$.

We next extend the previous integration by parts formula to matrix fields. Let $S$ and $T$
be two smooth matrix fields on $\mathbb{R}^3$ not necessarily symmetric.
We have
\begin{equation} \label{rel07}
\int_{\Omega} S: \curl T \ \ud x = \int_{\Omega} \curl S : T \ \ud x +
\int_{\partial \Omega} SC_n : T \ \ud x_f
\end{equation}
and
\begin{equation} \label{rel09}
\int_{f} \rotf S : T \ \ud x_f = \int_{f}S: \curlf T \ \ud x_f - \int_{\partial f}
(S s)\cdot(Tn) \ud s.
\end{equation}

Given a symmetric matrix field $S$ define $\Lambda_f(S):f\to Q_n\Sym Q_n$ by
\begin{equation*}
 \Lambda_f(S) = 2 \eps_f (Sn)- Q_n\partial_n S Q_n,
\end{equation*}
where $\partial_n S := \partial S / \partial n$. Hence, $\Lambda_f(S)$
is a symmetric matrix field defined on $f$.
If $T=\eps(v)$, where $v$ is a vector field then we have
\[
2 \eps_f (Tn) = \gradf \gradff(n'v) + Q_n\partial_n
\eps(v) Q_n.
\]
Hence, we obtain that
\begin{equation}\label{rel13}
\Lambda_f\bigl(\eps(v)\bigr) = \gradf \gradff(n'v).
\end{equation}
The tangential--normal components of the matrix field $\ccurl S$ on $f$
can be expressed in terms of $\Lambda_f(S)$.
Indeed, by the definition of the operator
$\rotf$ and
\eqref{rel06} we have
\begin{equation*}
C_n (\ccurl S)P_n = C_n \rotf \curll S
= \rotf C_n (\curlff S + \rotff S + C_n \partial_n S).
\end{equation*}
However, $C_n \rotff S = 0$ and, by \eqref{rel06}, $C_n \curlff S = Q_n
\gradff n'S$. Hence,
\begin{equation}\label{rot-Lambda}
C_n (\ccurl S)P_n = \rotf Q_n(\gradff n'S - \partial_n S)Q_n =
\rotf \Lambda_f(S).
\end{equation}
The next lemma indicates how the operator $\Lambda_f$ arises when integrating 
$\ccurl$ by parts.
\begin{lemma} \label{int-parts-curlcurl}
Let $S$ and $T$ be two smooth
matrix fields, with $S$ symmetric, on $\Omega$, where $\Omega$ is a bounded
subset of $\mathbb{R}^3$ with a
piecewise smooth boundary $\partial \Omega$. Then
\begin{multline*}
\int_{\Omega} S : \ccurl T \ud x=
 \int_{\Omega} \curll \curl S:T \ud x
- \int_{\partial \Omega} C_n S C_n: \frac{\partial T}{\partial n} \ud x_f\\
- \int_{\partial \Omega}[C_n \rotf S
- (\rotf^* S ) C_n
   + C_n \Lambda_f(S) C_n]: T \ud x_f.
\end{multline*}
\end{lemma}

\begin{proof}
Since $\partial \Omega$ has no boundary, it follows from \eqref{rel09} that
\begin{align}
\int_{\partial \Omega} \rotf S :T \ \ud x_f = \int_{\partial \Omega} S: \curlf T \ \ud x_f.
\end{align}
Then using \eqref{rel05} and the antisymmetry of $C_n$ we obtain
\begin{align} \label{rel19}
\int_{\partial \Omega} S: \curl T \ud x_f &= \int_{\partial \Omega} S:
\curlf T \ud x_f +
\int_{\partial \Omega} S: \rotf T  \ud x_f
- \int_{\partial \Omega} S: \frac{\partial T}{\partial n} C_n
\ud x_f\nonumber\\
& = \int_{\partial \Omega} (\rotf + \curlf)S:T \ud x_f +
\int_{\partial \Omega} S C_n :  \frac{\partial T}{\partial n} \ud x_f.
\end{align}
(So far we have not used the assumption that $S$ is symmetric.)
Next, from \eqref{rel07} we get
\begin{equation*}
\begin{split}
\int_{\Omega} S: \ccurl T \ud x &=
\int_{\Omega} \curl S : (\curl T')' \ud x
+ \int_{\partial \Omega} S C_n: (\curl T')' \ud x \\
& = \int_{\Omega} (\curl S)' : \curl T' \ud x
- \int_{\partial \Omega} C_n S': \curl T' \ud x_f \\
\end{split}
\end{equation*}
and
\[
\int_{\Omega} (\curl S)' : \curl T' \ud x = \int_{\Omega} \curl(\curl S)':T' \ud x +
\int_{\partial \Omega} (\curl S)' C_n : T' \ud x_f. \\
\]
Since $[\curl(\curl S)']' = \curll \curl S$
we obtain the identity
\begin{multline*}
\int_{\Omega} \curll \curl S: T \ud x 
\\
= \int_{\Omega} S: \ccurl  T \ud x
+ \int_{\partial \Omega} C_n \curl S : T \ud x_f
+ \int_{\partial \Omega} C_n S: \curl T' \ud x_f.
\end{multline*}
Using \eqref{rel19}, with $S$ replaced by $C_nS$,
we can rewrite the last term as
\begin{align*}
\int_{\partial \Omega} C_n S: \curl T' \ud x_f & =
\int_{\partial \Omega} C_n (\curlf + \rotf)S:T' \ud x_f +
\int_{\partial \Omega} C_n S C_n: \frac{\partial T'}{\partial n} \ud x_f \\
& = - \int_{\partial \Omega}[(\curlff+\rotff)S]C_n : T \ud x_f
+ \int_{\partial \Omega} C_n S C_n: \frac{\partial T}{\partial n} \ud x_f.
\end{align*}
Thus we obtain
\begin{multline}\label{rel999}
\int_{\Omega} S : \curl \curll T \ud x =
 \int_{\Omega} \curll \curl S:T \ud x
- \int_{\partial \Omega} C_n S C_n: \frac{\partial T}{\partial n} \ud
 x_f
\\
 - \int_{\partial \Omega} C_n \curl  S : T \ud x_f +
\int_{\partial \Omega}[(\curlff+\rotff)S]C_n : T \ud x_f.
\end{multline}
In order to see that \eqref{rel999} is equivalent to the desired
identity it is enough to show that
\begin{equation}\label{rel999b}
C_n \curl S - (\curlff S) C_n = C_n \Lambda_f(S) C_n + C_n \rotf S.
\end{equation}
However, from \eqref{rel7} it follows that
$C_n \curlf S - (\curlff S) C_n = 2 C_n \eps(Sn) C_n$,
and hence \eqref{rel999b} follows from \eqref{rel05}.
\end{proof}

The lemma can be used to determine when a
symmetric matrix field $S$ on $\Omega$ which is piecewise smooth
with respect to a given triangulation belongs to the space
$H(\ccurl, \Omega;\Sym)$. This holds if and only if $\ccurl S$ defined
piecewise coincides with $\ccurl S$ defined in the sense of distributions.
Integrating against a smooth test function and using the lemma, we
find that a necessary condition is that $C_n S C_n$, or equivalently $Q_n S Q_n$, is continuous
across element faces.  In this case the quantities
$C_n \rotf S$ and $(\rotff S)C_n$ are also continuous.  Again invoking
the lemma again, we obtain the following result.

\begin{theorem}\label{piecewise-smooth}
Assume that $S \in L^2(\Omega;\Sym)$ is piecewise smooth with respect
to a triangulation $\T$ of $\Omega$. Then $S \in H(\ccurl, \Omega;\Sym)$
if and only if $Q_n S Q_n$ and $\Lambda_f(S)$ are continuous across each
face in $\Delta_2(\T)$.
\end{theorem}

\section{Polynomial matrix fields on a single tetrahedron} \label{single-tet}
In order to study the space $\Sigma_h \subset H(\div ,\Omega;\Sym)$
introduced above, we will
need to study finite element subspaces $\Theta_h$ of $H( \ccurl,\Omega;\Sym)$.
The discussion in the present section will be
restricted to a single tetrahedron $K$, but based on the results
derived here we will define the space $\Theta_h$ in the final section of
the paper. In fact, we will present a complete discrete
elasticity complex of the form
\begin{equation*}
\mathbb{T} \hookrightarrow W_h \xrightarrow{\eps} \Theta_h
\xrightarrow{\ccurl}
 \Sigma_h
\xrightarrow{\div } V_h \to 0,
\end{equation*}
where $W_h \subset H^1(\Omega;\R^3)$ and $\Theta_h \subset
H(\ccurl,\Omega;\Sym)$
are piecewise polynomial spaces with respect to the triangulation $\T_h$.

Let $K \subset
\mathbb{R}^3$
be a fixed tetrahedron and define the polynomial space
\begin{equation}\label{defN}
\N_k = \N_k(K) = \{ \, S \in
\mathcal{P}_k(K;\mathbb{S}) \, | \, Q_n S Q_n =\Lambda_f(S)=0, \,
f \in \Delta_2(K)\}.
\end{equation}
Most of the discussion in this section is devoted to computing the
dimension and a basis for this space. However, first we need some
additional notation.

\begin{figure}[ht]
\centerline{%
\begin{picture}(0,0)%
\includegraphics{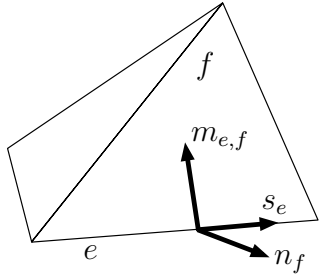}%
\end{picture}%
\setlength{\unitlength}{3947sp}%
\begin{picture}(1974,1725)(-11,-874)
\put(1151,-31){\makebox(0,0)[lb]{\smash{{\SetFigFont{12}{14.4}{\familydefault}{\mddefault}{\updefault}{\color[rgb]{0,0,0}$m_{e,f}$}%
}}}}
\put(1673,-805){\makebox(0,0)[lb]{\smash{{\SetFigFont{12}{14.4}{\familydefault}{\mddefault}{\updefault}{\color[rgb]{0,0,0}$n_f$}%
}}}}
\put(476,-770){\makebox(0,0)[lb]{\smash{{\SetFigFont{12}{14.4}{\familydefault}{\mddefault}{\updefault}{\color[rgb]{0,0,0}$e$}%
}}}}
\put(1182,391){\makebox(0,0)[lb]{\smash{{\SetFigFont{12}{14.4}{\familydefault}{\mddefault}{\updefault}{\color[rgb]{0,0,0}$f$}%
}}}}
\put(1593,-458){\makebox(0,0)[lb]{\smash{{\SetFigFont{12}{14.4}{\familydefault}{\mddefault}{\updefault}{\color[rgb]{0,0,0}$s_e$}%
}}}}
\end{picture}%
}
\caption{The $(n_f,s_e,m_{e,f})$ coordinate system for a face
$f$ and edge $e$ of the tetrahedron $K$.}
\label{fg:tet}
\end{figure}

If $f\in\Delta_2(K)$, we denote by $h_f$ the
perpendicular distance from the opposite vertex to $f$ and by $n=n_f$
the outward normal vector to $f$.   If $e$ is an edge, we let $s=s_e$
denote one of the unit vectors parallel to $e$.  When the edge $e$
belongs to the face $f$, we write $m=m_{e,f}$ for the unit vector in $f$,
normal to $e$, pointing from $e$ into $f$.  See Figure~\ref{fg:tet}.
When the notations $f_+$ and $f_{-}$ are used to denote two faces,
the corresponding normals will be denoted $n_+$ and $n_{-}$, and the
perpendicular distances $h_+$ and $h_{-}$, respectively.  The notations
$m_+$ and $m_{-}$ will also be used to denote $m_{e,f_+}$ and $m_{e,f_-}$
where $e$ is the edge common to $f_+$ and $f_-$.

The barycentric coordinates on $K$ will be labelled by the faces.  That is,
they are
$\lambda_f \in \mathcal{P}_1(K;\mathbb{R})$ determined by $\lambda_f \equiv 0$ on $f$ and
$\sum_f \lambda_f \equiv 1$  on $K$.
We recall that grad $\lambda_f=-n_f/h_f$.
Let $g \in \Delta(K)$ be a face of dimension $m$ with vertices
$x_{i_0},x_{i_1},\ldots,x_{i_m}$. For
$m=0$, $g$ is a vertex, for $m=1, g$ is an edge, and so on.
We define the bubble functions $b_g = \lambda_{f_{i_0}}
\lambda_{f_{i_1}}\cdots\lambda_{f_{i_m}}$,
where
$f_{i_k}$ is the face opposite vertex $x_{i_k}$. For
$d>0$,
$\mathcal{P}^g_d(K)=\operatorname{span} \{ \,
\lambda_{f_{i_0}}^{j_0}\lambda_{f_{i_1}}^{j_1} \cdots
\lambda_{f_{i_m}}^{j_m} \, | \, j_0+\cdots+j_m=d \, \}$,
so that $\dim \mathcal{P}^g_d(K) = \dim \mathcal{P}_d(g)$.
Note that if $x$ is a vertex opposite face $f$, then $b_x=\lambda_f$ and
$\mathcal{P}^x_d(K) =\R\lambda_f^d$,
while
$b_K=\prod_f \lambda_f$ and $\mathcal{P}^K_d = \mathcal{P}_d$. For
a given face $f_0 \in \Delta_2(K)$,
$b_{f_0} = b_K/\lambda_{f_0}=\prod_{f \neq f_0} \lambda_f$,
and for a given edge $e \in \Delta_1(K)$,
$b_e = b_K/(\lambda_{f_{-}} \lambda_{f_+} )$, where $f_-$ and $f_+$ are
the faces containing $e$.

The monomials of degree $k$ in the barycentric coordinates form a basis for
$\mathcal P_k(K)$, and by grouping together terms according to which coordinates
enter the monomial, we can uniquely represent any $p\in\mathcal P_k$ as
\begin{equation}\label{Lagrange}
p = \sum_{ g \in \Delta(K) } b_g p_g,\quad p_g \in \mathcal{P}^g_{k-1-\dim g}(K).
\end{equation}
The standard Lagrangian degrees of freedom for $p\in\mathcal P_k(K)$ are the values
of $p$ at the vertices, the moments of $p$ of degree at most $k-2$ on each of the
edges of $K$, the moments of degree at most $k-3$ on each of the faces, and
the moments of degree at most $k-4$ on $K$.  From the vertex values of $p$
we may determine the polynomials $p_g$ in \eqref{Lagrange}
for $g\in\Delta_0(K)$.  From these and the edge moments
we may determine as well the $p_g$ for $g\in\Delta_1(K)$, etc.

Of course analogous considerations apply to $\mathcal P_k(K;X)$ for
$X$ a vector space.  In particular, we have the representation
$p = \sum_{ g \in \Delta(K) } b_g p_g$ for $p\in \mathcal P_k(K;X)$
where now $p_g\in \mathcal{P}^g_{k-1-\dim g}(K;X)$, the space of
$X$-valued polynomials on $K$ whose components with respect to a
basis of $X$ belong to $\mathcal{P}^g_{k-1-\dim g}(K)$.

If $e = f_- \cap f_+ \in \Delta_1(K)$, with $f_-, f_+ \in
\Delta_2(K)$, we
let $G_e \in \Sym$ be the matrix
\[
G_e = n_- n_+' + n_+ n_-'.
\]
We note that $G_e s = 0$ and $m_-'G_em_- = m_+'G_em_+ = 0$.

\begin{lemma} \label{lemT1}
For $k \geq 0$, the dimension of the space
$$
\N_k^0 = \N_k^0(K) := \{ \, S \in \mathcal{P}_k(K;\mathbb{S}) \, | \, Q_n S Q_n =0,\quad
f\in\Delta_2(K)\,\}
$$
is $(k+1)k(k-1)$.
\end{lemma}
\begin{proof}
For a subsimplex $g \in \Delta(K)$, let us first define
$$
N^g =\{\, S \in \Sym \, | \, Q_n S Q_n = 0 
\text{ for all faces containing $g$}  \,\}.
$$
Clearly, if $g=K$ then dim $N^g=6$, if $g \in \Delta_2(K)$ then
dim $N^g=3$, and if $g \in \Delta_0(K)$ then dim $N^g=0$.
Finally, if $g \in \Delta_1(K)$, i.e., $g =e$ is an edge, then dim
$N^g=1$.
In fact, the space $N^e$ is then spanned by the matrix $G_e$ introduced above.

If $S \in \N_k^0$ then from \eqref{Lagrange}
we obtain the representation
\begin{equation} \label{N0-rep}
S = \sum_{g \in \Delta(K)} b_g S_g, \quad S_g \in
\mathcal{P}^g_{k-1-\dim g}(K;N^g).
\end{equation}
As a consequence
\begin{align*}
\dim \N^0_k &= \sum_{g \in \Delta(K)}
\dim \mathcal{P}^g_{k-1-\dim g}(K;N^g) \\
& = 6(k-1) + 6(k-1)(k-2) + (k-1)(k-2)(k-3)  = (k+1)k(k-1).
\qedhere
\end{align*}
\end{proof}
Before we are able to compute $\dim\N_k$, we need to establish several lemmas.
\begin{lemma} \label{zeroedge}
If $S \in \N_k$, $S$ is zero on each edge.
\end{lemma}
\begin{proof}
As above let $e = f_- \cap f_+ \in \Delta_1(K)$, with $f_-, f_+ \in
\Delta_2(K)$. Then
\[
m_-'n_+ = m_+' n_- < 0.
\]
In fact, the transformation
\begin{equation}\label{rotation}
\begin{pmatrix} m_-\\ n_- \end{pmatrix}
\mapsto
\begin{pmatrix} n_+\\ m_+ \end{pmatrix}
\end{equation}
is a rotation in the plane orthogonal to $e$.

Let $S \in \N_k$.
Since $\N_k \subset \N_k^0$
we know that $S|_e = \rho G_e$, where $\rho \in \P_k(e)$.
Therefore, if we can show that
\begin{equation}\label{desire}
m_+'Sn_+ + m_-'Sn_- = 0
\end{equation}
on $e$, then $\rho (m_+'n_- + m_-'n_+) = 2 \rho m_+'n_- = 0$,
and, as a consequence, $S$ is zero on $e$. It therefore suffices to
show \eqref{desire}.

On $e$, we must have $s' \Lambda_f(S) m = 0$, i.e.,
$$
\partial_s(m'Sn) + \partial_m(s'Sn) - \partial_n(s'Sm)=0,
$$
where $f$ is either $f_-$ or $f_+$. By adding this property for
the two faces we obtain
\begin{multline*}
- \partial_s( m_{+}' S n_+ + m_{-}' S n_{-} )
\\
= [ \partial_{m_{-}}(s'Sn_{-}) - \partial_{n_{-}}(s'Sm_{-})]
 - [\partial_{n_+}(s'Sm_+) -\partial_{m_{+}}(s'Sn_{+}) ].
\end{multline*}
However, the right hand side here is zero as a consequence of the fact
that the transformation \eqref{rotation} is a rotation. In fact, this property
implies that
\[
\partial_{m_{-}}(v'n_{-}) - \partial_{n_{-}}(v'm_{-})
=
\partial_{n_+}(v'm_+) -\partial_{m_{+}}(v'n_{+})
\]
for any smooth vector field $v$ on $K$. Hence, we can conclude that
$m_{+}' S n_+ + m_{-}' S n_{-}$ is a constant along $e$, and
since it is zero at the vertices, \eqref{desire} holds.
\end{proof}

\begin{lemma}\label{dimnkK}
Let\begin{multline*}
\N_{k,\partial K} :=\{\, U \in \N_k^0 \, | \,
U=\sum_{f \in \Delta_2(K)} b_f U_f ,\, U_f \in \P_{k-3}^f(K; \Sym) \text{ and}
\\ \Lambda_f(U)|_{e} =0 \text{ for each face $f$ and each edge $e$ of $f$}\, \}.
\end{multline*}
Then $\dim\N_{k,\partial K}= 6(k^2-6k+10)$.
\end{lemma}
\begin{proof}
If $U = \sum_{f} b_f U_f$, then $U \in \N_k^0$ if and only if
each coefficient $U_f \in \P_{k-3}^f(K; \Sym)$ satisfies
\begin{equation*}
Q_{n} U_f Q_{n} = 0\text{\quad on $f$}.
\end{equation*}
Hence, this property is assumed to hold.
We have $\Lambda_f(U)=0$ on an edge $e\subset f$ if and only if the
three terms
$s'\Lambda_f(U)s$, $s'\Lambda_f(U)m$ and $m'\Lambda_f(U)m$ vanish there.
For any fixed unit vector $t$ and $e \in \Delta_1(K)$, we have
\begin{equation}\label{edgrel}
\partial_t U = \sum_f (\partial_t b_f)U_f = -b_e (\frac{t' n_+}{h_+} U_{f_{-}}
+ \frac{t'n_{-}}{h_{-}} U_{f_{+}}) \text{\quad on $e$},
\end{equation}
where $f_-$ and $f_+$ are the two faces meeting the edge $e$
and we have used that $\grad b_{f_{-}} = -n_+b_e/h_+$ on $e$.

Recall that $ s'\Lambda_f(U)s = 2 \partial_s(s'Un) - \partial_n (s'Us) $.
Since $U=0$ on $e$,
$$
s'\Lambda_f(U)s  = - \partial_n (s'Us) \text{\quad on $e$}.
$$
However, since $Q_{n} U_f Q_{n}=0$ on the face $f$, we have $ s'
U_{f_{-}} s =s'U_{f_{+}}s=0$ on $e$. By \eqref{edgrel}, with $t=n$,
we conclude that $s'\Lambda_f(U)s=0$ on $e$.

Next, similar considerations for $ s'\Lambda_f(U)m =\partial_s(m'Un)
+ \partial_m (s'Un)-\partial_n (s'Um)$ give
\begin{equation*}
s'\Lambda_{f_+}(U)m_+  = \partial_{m_+}(s'Un_+) -
\partial_{n_+} (s'Um_+)  \text{\quad on $e$}.
\end{equation*}
Furthermore, from \eqref{edgrel} we have
\[
\partial_{m_+}(s'Un_+) = - m'_+n_{-} b_e \frac{s' U_{f_{+}}n_+}{h_{-}}
\]
on $e$, and, using $s'U_{f_{-}}m_-=s'U_{f_+}m_+ = 0$, we obtain
\[
\partial_{n_+} (s'Um_+) =
-b_e \frac{s'U_{f_{-}}m_+}{h_{+}}
=-m'_+ n_{-}b_e\frac{s'U_{f_{-}}n_-}{h_{+}}.
\]
It follows that
$$
s'\Lambda_{f_+}(U)m_+  = m'_+n_{-} b_e ( \frac{s'U_{f_{-}}n_{-}}{h_+} -
\frac{s' U_{f_{+}}n_+}{h_{-}}) \text{\quad on $e$}.
$$
We have therefore shown that
$s'\Lambda_{f}(U)m=0$ on all edges if and only if,
\begin{equation}\label{lkb1}
\frac{s'U_{f_{-}}n_{-}}{h_+} = \frac{s' U_{f_{+}}n_+}{h_{-}}
\end{equation}
on all edges of $K$.

Finally, we consider $m'\Lambda_f(U)m = 2 \partial_m (m'Un)
-\partial_n (m'Um)$.
Using the fact that both $m_-'U_{f_-}m_-$ and $m_+'U_{f_+}m_+$ vanish on
$e$,
we obtain from \eqref{edgrel} that, on $e$,
\begin{align*}
\partial_{n_+} m_{+}' U m_{+} &= -b_e
\frac{m_{+}' U_{f_{-}} m_+}{h_{+} }
\\
&=-\frac{b_e}{h_+} [(m'_+ n_{-})^2 n_{-}'U_{f_{-}}n_{-} +2m'_+ m_{-}m'_+ n_{-}
m_{-}'U_{f_{-}}n_{-}],
\end{align*}
so
\begin{equation*}
m_{+}' \Lambda_{f_+}(U) m_+
= m'_+ n_{-}b_e[ m'_+ n_{-} \frac{n_{-}'U_{f_{-}}n_{-}}{h_+} -2
( \frac{m_{+}' U_{f_+} n_+}{h_{-}} - m'_+ m_{-}
\frac{m_{-}'U_{f_{-}}n_{-}}{h_+})].
\end{equation*}
Hence,  $m_{+}' \Lambda_{f_+}(U) m_+$ vanishes on $e$
if and only if,
\begin{equation*}
n_{-}'U_{f_{-}}n_{-} = \frac{2 h_+}{m'_+ n_{-}} (\frac{m_{+}' U_{f_+}n_+}{h_{-}} -
m'_+ m_{-} \frac{m_{-}'U_{f_{-}}n_{-}}{h_+} ) \text{\quad on $e$}.
\end{equation*}
Note that this condition is not symmetric in $f_-$ and $f_+$.
We thus obtain two conditions for each edge $e$.

Since $U_f \in \P_{k-3}^f(K; \Sym)$, with $ Q_{n} U_f Q_{n}=0$
on $f$, it follows that $U_f$ is uniquely determined by
the vector field $v_f := U_f n \in \P_{k-3}(f; \R^3)$.
The analysis above shows that
$U = \sum_f b_f U_f \in \N_{k,\partial K} $
if and only if these vector fields satisfy
\begin{itemize}
\item[(A)] $\displaystyle \frac{s' v_{f_{-}} }{h_+} = \frac{s'v_{f_+}}{h_{-}}$ on $e$, and
\item[(B)] $\displaystyle n'_{-} v_{f_{-}}  =
\frac{2 h_+}{m'_+ n_{-}} (\frac{m_{+}' v_{f_+}}{h_{-}} -
m'_+ m_{-} \frac{m_{-}' v_{f_{-}} }{h_+})$ on $e$,
\end{itemize}
whenever an edge $e$ is shared by faces $f_-$ and $f_+$.
Therefore, there is an isomorphism between $\N_{k,\partial K}$ and
\begin{equation}\label{isomorphic}
\{\,(v_f) \in \prod_{f \in \Delta_2(K)}
\mathcal{P}_{k-3}(f;\mathbb{R}^3)\,|\, \text{the $v_f$ satisfy (A) and (B)}\,\}.
\end{equation}
To compute the dimension of the space
\eqref{isomorphic} we consider the relations (A) and (B)
at a fixed vertex $x$ of $K$. If $(v_f)$ is an element of the space
\eqref{isomorphic} define $z \in \R^3$ by
\begin{equation}\label{tangential-comp}
s'z = h_f s'v_f(x)
\end{equation}
for $s$ chosen as tangents to each edge $e$ meeting $x$, and where
$f$ is a face meeting $e$. Note that the vector $z$ is well--defined as
a consequence of condition (A), and that
for each face $f$ containing $x$ we have
$h_f t' v_f(x) = t'z$
for all vectors $t$ which are tangential to the face $f$. 
Using the expansion $n_-=[m_+-(m'_+ m_)-m_-]/m'_+ n_-$ we
can then rewrite condition (B) at the vertex $x$ as
\begin{equation}\label{normal-comp}
h_f n'_fv_f(x)  = 2 n'_fz.
\end{equation}
From this discussion we can conclude that the dimension
of the space \eqref{isomorphic}, and hence
$\dim\N_{k,\partial K}$, is at least $6(k^2 - 6k + 10)$. To see
this observe that
$$
\dim \prod_{f}
\mathcal{P}_{k-3}(f,\mathbb{R}^3) = 6(k-2)(k-1).
$$
Furthermore, the
conditions (A) and (B) represent a total of $6\cdot 3 \cdot
(k-4)
= 18(k-4)$ constraints in the interior of the edges and $4 \cdot 6 =
24$
constraints at the vertices. Since
$6(k-2)(k-1) - 18(k-4) - 24 = 6(k^2 - 6k + 10)$,
this is a lower bound for $\dim\N_{k,\partial K}$.

We complete the proof by showing
that elements of the space
\eqref{isomorphic} are determined by $6(k^2 - 6k + 10)$
degrees of freedom, in fact
by degrees of freedom corresponding to the space
\[
\prod_{x \in \Delta_0(K)} \R^3
\times
\prod_{e \in \Delta_1(K)} \mathcal{P}_{k-5}(e;\mathbb{R}^3) \times
\prod_{f \in \Delta_2(K)} \mathcal{P}_{k-6}(f;\mathbb{R}^3).
\]
To see this,
for each vertex $x$ pick a vector $z=z(x) \in \R^3$
and choose $v_f(x)$ such that
the relations \eqref{tangential-comp} and \eqref{normal-comp} hold
for all faces meeting $x$. This determines the vectors $v_f(x)$
for all vertices $x \in \partial f$.
We then define $h_f s'v_f $ on each edge by the standard interior
degrees of freedom, and $m'v_f$ with respect to
both faces meeting $e$ are determined similarly.
These degrees of freedom on the edges correspond to the space
$\prod_{e} \mathcal{P}_{k-5}(e;\mathbb{R}^3)$. The normal components
$n_f'v_f$ are determined on each edge by
condition (B). Finally, we must apply the interior degrees of
freedom to $v_f$ on $f$.
We conclude that elements of the space $ \N_{k,\partial K} $
are determined by
$12 + 18\dim\mathcal{P}_{k-5}(e) + 12 \dim
\mathcal{P}_{k-6}(f) = 6(k^2-6k+10)$
degrees of freedom.
\end{proof}

For $U\in\N_{k,\partial K}$ and $f\in\Delta_2(K)$,  $\Lambda_f(U)$ is a polynomial vanishing on $\partial f$,
and so the quotient $\Lambda_f(U)/b_f$ is a polynomial. We define
$T_f:\N_{k,\partial K}\to \mathcal P_{k-4}(f,Q_n\Sym Q_n)$
by
$$T_f(U)=-h_f \Lambda_f(U)/b_f.
$$
\begin{lemma}\label{edgecompa}
If $f_-,f_+\in\Delta_2(K)$, $e=f_-\cap f_+$, and $s$ is a unit vector parallel to $e$, then
\begin{equation*}
s' T_{f_+}(U) s = s' T_{f_{-}}(U) s \text{\quad on $e$},\quad
U\in\N_{k,\partial K}.
\end{equation*}
\end{lemma}

\begin{proof}
First we show that
\begin{equation}\label{dmsls}
\partial_{m_{+}} s'\Lambda_{f_{+}}(U) s = \partial_{m_{-}} s'\Lambda_{f_{-}}(U) s
\text{\quad on $e$}.
\end{equation}
Recall that $\partial_{m_+} b_{f_+} = - m'_+ n_{-}b_e/h_{-}$,
$\partial_{m_{-}}
b_{f_{-}} = - m'_{-}n_{+}b_e/h_+$ on $e$ and $m'_+n_{-} = m'_{-}n_{+}$.
We have on an edge $e$,
$\partial_m s' \Lambda_f(U) s = 2 \partial_s \partial_m s' U n - \partial_m \partial_n s' U s$.
Using \eqref{edgrel} we obtain
$$
\partial_{m_{+}} s'U n_{+} = - b_e \frac{m'_+ n_{-}}{h_{-}} s' U_{f_+} n_+,
$$
which is symmetric in $f_{-}$ and $f_{+}$ as a consequence of
\eqref{lkb1} and $m'_+ n_{-} = m'_{-}n_{+}$. 

The identity \eqref{dmsls} will follow if we show that 
$\partial_{m_{+}} \partial_{n_{+}} s'U s
= \partial_{m_{-}} \partial_{n_{-}} s'U s$.
Consider first the term
\[
V= b_{f_-}U_{f_-} + b_{f_+}U_{f_+}.
\]
Since $Q_n U_f Q_n = 0$ for $f=f_-,f_+$ and 
$\grad b_{f_{-}} = -\frac{b_e}{h_+}n_+$ on $f_+$
we derive that at the edge $e$
\begin{align*}
\partial_{m_{+}} \partial_{n_{+}} s'V s 
&= -b_e[\frac{1}{h_+}s'\partial_{m_{+}}U_{f_-}s
+ \frac{m_+'n_-}{h_-}s'\partial_{n_{+}}U_{f_+}s]\\
&= -  (m_+'n_-)b_e[\frac{1}{h_+}s'\partial_{n_{-}}U_{f_-}s
+ \frac{1}{h_-}s'\partial_{n_{+}}U_{f_+}s],
\end{align*}
and this expression is symmetric in $f_-$ and $f_+$.
Finally, consider terms of the form 
$W = b_f U_f$, 
where $f$ is neither $f_-$ nor $f_+$. In this case 
$b_f = \lambda_{f_-} \lambda_{f_+} \lambda$, where $\lambda$ is the
barycentic coordinate associated the fourth face ($\neq f,f_-,f_+$) of
$K$,
and 
on $e$ we have
\[
\partial_{m_{+}} \partial_{n_{+}} s'W s 
= \frac{m_+'n_-}{h_-h_+}\lambda s'U_fs. 
\]
This is again symmetric in $f_-$ and $f_+$.
We have therefore established
\eqref{dmsls}.

Now, by definition, $h_+\Lambda_{f_+}(U)=-b_{f_+}T_{f_+}(U)$.  Therefore
\begin{equation*}
(m'_+ n_-)^{-1}h_+ h_- \partial_{m_+}s'\Lambda_{f_+}(U)s
= b_e s'T_{f_+}(U)s\text{\quad on $e$}.
\end{equation*}
By \eqref{dmsls}, the left-hand side is unchanged if we interchange the subscripts
$+$ and $-$, so the same must be true of the right-hand side.
\end{proof}

\begin{lemma}\label{techlemma}
Let $(T_f) \in \prod_{f \in \Delta_2(K)} \P_k(f; Q_n \Sym Q_n)$
be such that $s'T_{f_-}s=s'T_{f_+}s$ on $e$, whenever $e=f_-\cap f_+$,
$f_-,f_+\in\Delta_2(K)$.
Then there exist an $S \in \P_k(K; \Sym)$ such that $Q_n S Q_n = T_f$
for all $f \in \Delta_2(K)$.
\end{lemma}

\begin{proof}
We will define $S \in \P_k(K; \Sym)$ by first specifying its
vertex values, then specifying its moments of degree at most
$k-2$ on the edges, then its moments of degree at most $k-3$ on
faces, and then the moments of degree at most $k-4$ over the interior of $K$.

Let $x$ be a vertex.  We define the matrix $S(x)\in\Sym$ by specifying
the values $s'_iS(x)s_j$ where the $s_i$ are the tangents
to the edges $e_i$ meeting at $x$ (and so the $s_i$ form a basis for
$\R^3$).  Namely we take $s'_i S(x) s_j = s'_i T_f s_j$ with $f$
the face containing $e_i$ and $e_j$.  If $i=j$ there are two
possible choices of the face $f$, but they give the same result by assumption.

For the interior degrees of freedom on an edge $e$ we use the
basis $s$, $m_-$, $m_+$ of $\R^3$, and let
$T_e \in \P_k(e; \Sym)$ be given by
\begin{align*}
s'T_es = s'T_{f_-}s = s'T_{f_+}s, \quad s'T_em_- &= s'T_{f_-}m_-,
\quad s'T_em_+ = s'T_{f_+}m_+,\\
m_-'T_em_- = m_-'T_{f_-}m_-, \quad m_+'T_em_+ &= m_+'T_{f_+}m_+,
\quad m'_-T_e m_+ = 0.
\end{align*}
Then we define $S|_e$ by
\[
\int_e (S-T_e)V \, ds = 0, \quad V \in \P_{k-2}(e; \Sym).
\]
Similarly, for the interior degrees of freedom on each face
we let $Q_nSQ_n$ inherit the moments from $T_f$, while the data for
$SP_n$
is taken to be zero.
The interior
degrees of freedom on $K$ are all taken to be zero.
\end{proof}

As a consequence of the two previous lemmas, there is a map
\begin{equation}\label{mapS}
\N_{k,\partial K}  \to \P_{k-4}(K;\Sym), \quad U \mapsto
S(U),
\end{equation}
such that $Q_n S(U) Q_n = T_f(U)$ for all faces $f \in \Delta_2(K)$.

We are finally ready to compute the dimension of the space $\N_k$ defined
in \eqref{defN}.

\begin{theorem}\label{dim-N}
For $k \geq 3$, the dimension of the space
$\N_k$ is $k(k^2-6k+11)$.
\end{theorem}
\begin{proof}
 Let $S \in \N_k$. By Lemma~\ref{zeroedge}, $S$ must be zero on each edge and so
can be written $S= \sum_f b_f S_f + b_K S_K$,
where $S_f \in \mathcal{P}_{k-3}^f(K;\mathbb{S}) $ and $S_K \in \mathcal{P}_{k-4}(K;\mathbb{S})$.
Now $\eps_f(b_KS_Kn_f)$ vanishes on $f$ since $b_K$ does, while
$$
\partial_n(b_KQ_nS_KQ_n)=(\partial_n b_K)Q_nS_KQ_n=-\frac{b_f}{h_f}Q_nS_KQ_n \text{\quad on $f$}.
$$
Thus
\begin{equation}\label{Lambda_f-comp}
\Lambda_f(b_K S_K) =\frac{ b_f}{h_f} Q_n S_K Q_n.
\end{equation}
In particular, $\Lambda_f(b_K S_K)$ vanishes on $\partial f$.
It follows
that if $S \in \N_k$ and we define $U=\sum_f b_f S_f$
then $U \in \N_{k,\partial K}$.
Therefore, the map $(U,S_K)\mapsto U+b_KS_K$ defines an isomorphism
from
\begin{equation*}
\{\, (U,S_K) \in \N_{k,\partial K} \times \mathcal{P}_{k-4}(K;\Sym) \, | \, Q_{n} S_K Q_{n} = T_f(U)\text{\quad on each face $f$}\,\}
\end{equation*}
onto  $\N_k$.

Finally, note that a matrix field of the form $b_K V$, $V
\in \P_{k-4}(K;\Sym)$, belongs to $\N_k$ if and only if
$V \in \N^0_{k-4}$.
Therefore, using the map \eqref{mapS}, the mapping
\begin{equation*}
\N_{k,\partial K} \times \N^0_{k-4} \to
 \N_k, 
\quad
(U,V) \mapsto (U,S(U)+V),
\end{equation*}
is an isomorphism.
It follows that dim $\N_k=$ dim $\N_{k,\partial K}$ + dim $\N^0_{k-4}$
and using Lemma~\ref{dimnkK} and Lemma~\ref{lemT1} we get
$\dim \N_{k} = 6(k^2-6k+10) + (k-3)(k-4)(k-5) = k(k^2-6k+11)$.
\end{proof}

\section{The space of divergence-free matrix fields with vanishing normal traces}
\label{basis}
Recall that the space
\[
\M_k = \M_k(K) = \{\, S \in
\mathcal{P}_k(K;\mathbb{S}) \, |
\, \div S =0 \text{\quad on $K$},\quad
P_n S = 0, \quad f\in \Delta_2(K)\, \}
\]
appears in the degrees of freedom for the finite element space
$\Sigma_h \subset H(\div, \Omega; \Sym)$ introduced in Sections~\ref{lowest-order}
and \ref{higher-order}. Therefore, a derivation of the dimension of this space
is fundamental for our theory, while a construction of a (dual) basis
for the space $\M_k$ is necessary for the implementation of
the method.  The dimension formula will be a simple consequence
of the following lemma, in which $\P^0_{k+3}(K;\R^3) := \{ v \in \P_{k+3}(K;\R^3) \, | \,
v \equiv 0 \text{\quad on $\partial K$}\, \} = b_K \P_{k-1}(K;\R^3)$.

\begin{lemma} \label{problemma}
\begin{enumerate}
\item The operator $\ccurl$ maps $\N_{k+2}(K)$ onto
$\M_k(K)$.
\item $\{ \, T \in \N_{k+2} \, | \, \ccurl T = 0 \, \}
= \eps[\P^0_{k+3}(K;\R^3)]$.
\item The following sequence is exact:
\begin{equation}\label{ses}
0\to \P^0_{k+3}(K;\R^3)\xrightarrow{\eps}\N_{k+2}(K)
\xrightarrow{\ccurl}\M_k(K)\to0.
\end{equation}
\end{enumerate}
\end{lemma}
\begin{proof}
It follows directly from \eqref{rel12} and \eqref{rot-Lambda}
that $\ccurl \N_{k+2} \subset  \M_k$.
Hence, to prove the first statement we need only show that $\M_k
\subset \ccurl \N_{k+2}$.
Let $S \in \M_k$. Since $\div S=0$, it follows from the
exactness
of the complex \eqref{cd1}  that there is a $T
\in \mathcal{P}_{k+2}(K; \mathbb{S})$ such that $S = \ccurl T$.
The proof will be completed by constructing a vector field
$u \in \mathcal{P}_{k+3}(K;\mathbb{R}^3)$ such that
\begin{equation}\label{desired2}
Q_n \bigl(T - \eps (u)\bigr)Q_n =0,\quad
\Lambda_f\bigl(T- \eps (u)\bigr) =0, \text{\quad on each face $f$}.
\end{equation}
Note that since $S \in \M_k$ it follows from \eqref{rel12} that
$$
\rotf \rotff Q_n T Q_n =P_n (\ccurl T) P_n = P_n S P_n = 0 \text{\quad on each face $f$}.
$$
Hence, from the exact sequence \eqref{pol-elas-2d1} we conclude that
for each face $f \in \Delta_2(K)$ there is a vector field
$v_f \in \P_{k+3}(f; \R^3)$, with $P_nv_f =0$,
such that $Q_n T Q_n = \eps_f(v_f)$.
The vector fields $v_f$ are
uniquely determined up to a 2D rigid motion and hence we
may normalize them so that $\int_e s'v_f\,d s=0$ on each edge
$e\subset f$.
Since
$\partial_s (s' v_{v_{f_{-}}}) = s'Ts = \partial_s (s' v_{f_+})$
on each edge, we obtain that
$P_s v_{f_{-}} = P_s v_{f_{+}}$
on each edge $e=f_+ \cap f_{-}$. As a consequence, there is
$v \in \mathcal{P}_{k+3}(K;\mathbb{R}^3)$ such that $Q_nvQ_n=v_f$ on
each face $f$.  Then
$Q_n\eps(v)Q_n=\eps_f(v_f)=Q_nTQ_n$, i.e.,
$Q_n U Q_n = 0$
on each face $f$, where $U = T - \eps (v)$.
This implies, in particular, that
$Us$ and $\grad(s' U s)$ vanish on each edge $e \in \Delta_1(K)$.
Therefore,
\begin{equation}\label{slt}
s' \Lambda_f(U) s = 2\partial_s(n'Us) - \partial_n(s'Us) = 0\text{\quad on $\partial f$}
\end{equation}
for each face $f$.

Next, observe that, by \eqref{rot-Lambda},
$\rotf \Lambda_f(U) = C_n S P_n = 0$.
Hence, \eqref{pol-elas-2d2} implies that there is a scalar field $q_f \in
\mathcal{P}_{k+4}(f;\mathbb{R})$, uniquely determined up to
a linear function on $f$, such that
$\Lambda_f(U) = \gradf \gradff q_f$.
On each edge $e \in \partial f$,
we have by \eqref{slt} that
$0 = s' \Lambda_f(U) s = \partial_s^2 q_f$.
It follows that we can assume that $q_f \equiv 0$ on $\partial f$.
Hence, there exists another
vector field $w$ in $ \mathcal{P}_{k+3}(K;\mathbb{R}^3)$ such that
$Q_n w = 0$ and $n'w = q_f$
on each face. Recall by \eqref{rel13} that $\Lambda_f\bigl(\eps (w)\bigr)
= \gradf \gradff q_f = \Lambda_f(U)$.
Hence, if we let $u = v + w$, then the relation \eqref{desired2} holds.
This proves the first statement.

We now prove the second statement.
If $T=\eps(v)$ for some $v\in\mathcal P^0_{k+3}(K;\R^3)$, then
$\ccurl T=0$, and, by \eqref{rel13},
\begin{equation}\label{tv}
Q_n T Q_n = \eps_f(Q_n v), \quad \Lambda_f(T) = \gradf \gradff (n'v)\text{\quad on $f$},
\end{equation}
for each face $f$.  Since $v$ vanishes on $f$, the right hand sides of
these equations vanish, and so $T$ belongs to $\N_{k+2}(K)$.
Conversely, if $T\in\N_{k+2}(K)$, then, by the exactness
of the sequence \eqref{cd1}, $T=\eps(v)$
for some $v\in\mathcal P_{k+3}(K;\R^3)$, which is determined uniquely
if we require that
\begin{equation}\label{gauge}
\int_e s'v \ud s = 0, \quad e \in \Delta_1(K).
\end{equation}
(The functionals $v\mapsto\int_e s'v\ud s$, $e\in\Delta_1(K)$, form a set
of degrees of freedom for $\mathbb T$, the null space of $\eps$.)
From the first equation in \eqref{tv} and \eqref{gauge}, we find that $Q_nv$
vanishes on each face $f$.  Therefore the entire vector $v$ vanishes on
each edge $e$.  Using the second equation in \eqref{tv}, we see that
$n'v$ vanishes on each face as well, so $v\in\P^0_{k+3}(K;\R^3)$.
This completes the proof of the second statement.

The third statement is an immediate consequence of the first two
and the fact that $\mathbb T\cap \P^0_{k+3}(K;\R^3)=0$.
\end{proof}

\begin{theorem} \label{dimension-formula}
For $k \geq 4$ the space $\M_k(K)$
has dimension
$ (k+2)(k-2)(k-3)/2$.
\end{theorem}
\begin{proof}
Using first the short exact sequence in the lemma and then the dimension
formula in Theorem~\ref{dim-N}, we get
 \begin{align*}
\dim \M_k  &= \dim \N_{k+2} -
\dim \eps[\P^0_{k+3}(K;\R^3)]
\\
& = (k+2)(k^2-2k+3) - \dim  \P_{k-1}(K;\R^3) = (k+2)(k-2)(k-3)/2.
\qedhere
\end{align*}
\end{proof}

To conclude this section, we construct a basis for  the space
$\M_k(K)$ for $k=4$ and $k=5$.
(Alternatively a basis could be constructed for any $k$
using computational algebra software.)  For this we use
the following lemma, similar to Lemma~\ref{problemma}.

\begin{lemma}\label{problemma2}
\begin{enumerate}
\item The operator $\ccurl$ maps $b_K\N^0_{k-2}(K)$ onto $\M_k(K)$.
\item $\{ \, T \in b_K\N^0_{k-2} \, | \, \ccurl T = 0 \, \}
= \eps[b_K^2\P_{k-5}(K;\R^3)]$.
\item The following sequence is exact:
\begin{equation}\label{ses2}
0\to b_K^2\P_{k-5}(K;\R^3)\xrightarrow{\eps}b_K\N^0_{k-2}(K)
\xrightarrow{\ccurl}\M_k(K)\to0.
\end{equation}
\end{enumerate}
\end{lemma}
\begin{proof}
Note that $b_K\N^0_{k-2}\subset\N_{k+2}$ by 
\eqref{Lambda_f-comp}, and so $\ccurl b_K\N^0_{k-2}
\subset \M_k(K)$.

First we prove 2.  Let $w\in\P_{k-5}(K;\R^3)$.  By the Leibniz rule
\begin{equation*}
\eps(b_K^2w) = b_K^2\eps(w)+b_K[(\grad b_K)w'+w(\grad b_K)'].
\end{equation*}
Clearly $b_K\eps(w)\in \N^0_{k-2}$, and,
recalling that $\grad b_K=-b_f n_f/h_f$, we see that
$(\grad b_K)w'+w(\grad b_K)'\in \N^0_{k-2}$.  Thus
$\eps(b_K^2w)\in b_K\N^0_{k-2}$, giving the inclusion $\supset$.
Conversely, if $T\in b_k\N^0_{k-2}$ with $\ccurl T=0$, then,
by Lemma~\ref{problemma}, $T=\eps(b_Kv)$ for some $v\in \P_{k-1}(K;\R^3)$
and we need to show that $v=0$ on $\partial K$. Using the Leibniz rule
and the fact that $T$ vanishes on $\partial K$, we get that $nv'+vn'$
vanishes on each face.  We conclude that $v$ vanishes on the face,
using the elementary identity $v= (I+Q_n)(nv'+vn')n/2$.

It follows that
\begin{equation*}
\dim[\ccurl( b_K\N^0_{k-2})]=\dim\N^0_{k-2}-\dim\P_{k-5}(K;\R^3)=
\dim\M_k,
\end{equation*}
where we have used Lemma~\ref{lemT1} and Theorem~\ref{dimension-formula}.
The exactness of \eqref{ses2}, and so also the first statement of the lemma,
follows.
\end{proof}

Thus for $k=4$, $\ccurl$ is injective on
$b_K \N^0_2$, and so a basis for
$\M_4 = \ccurl (b_K \N^0_{2})$ is computable
directly from a basis for $\N^0_{2}$, which may be
obtained directly from \eqref{N0-rep}. 

Now let $k=5$.  The map
$\ccurl$ is not injective on $b_K
\N^0_3$, but has a kernel of dimension $3$. In this case, the representation \eqref{N0-rep}
presents an arbitrary element $S\in\N^0_3$ as
\begin{equation*}
S=\sum_{e\in\Delta_1(K)} b_e S_e + \sum_{f\in\Delta_2(K)} b_f S_f, \quad S_e\in\P^e_1(K;N^e),
\ S_f\in N^f.
\end{equation*}
Fix a particular face $f_0\in\Delta_2(K)$ and define $\N^{00}_3$ as the subspace of
$S\in\N^0_3$ for which $S_{f_0}=0$ in this representation, clearly
a subspace of codimension $3$.  We claim that
$\ccurl$ is injective on the space $b_K \N^{00}_{3}$,
and hence a basis for $\M_5$
can be computed from a corresponding basis of $\N^{00}_{3}$.
The injectivity follows since if $w \in\R^3$, with $\eps(b_K^2 w) \in b_K
\N^{00}_{3}$, then we get $w = 0$ arguing as in the proof of
Lemma~\ref{problemma2}.

\section{A Discrete Elasticity Complex}
\label{elasticity-complex}
The results of Section~\ref{basis} above completes the description of
the finite element spaces $V_h$ and $\Sigma_h$, introduced in Sections~\ref{lowest-order}
and \ref{higher-order}, and therefore
also of the finite element method \eqref{system}.
However, as already indicated in the beginning of Section~\ref{single-tet},
there are more structures hidden in the construction above. In fact,
the spaces $V_h$ and $\Sigma_h$ are constituents of a discrete
elasticity complex of the form
\begin{equation} \label{cd-discrete}
\mathbb{T} \hookrightarrow W_h \xrightarrow{\eps} \Theta_h
\xrightarrow{\ccurl}
 \Sigma_h
\xrightarrow{\div } V_h \to 0,
\end{equation}
where $W_h \subset H^1(\Omega;\R^3)$ and $\Theta_h \subset
H(\ccurl,\Omega;\Sym)$
are piecewise polynomial spaces with respect to the triangulation $\T_h$.
Furthermore, there exist interpolation operators $\Pi_h^W$,
$\Pi_h^\Theta$, $\Pi_h^\Sigma$ and $\Pi_h^V$ such that the diagram
\begin{equation} \label{cd-final}
\begin{CD}
\mathbb{T} \hookrightarrow \ @. C^{\infty}(\Omega;\mathbb{R}^3) @>\eps >> C^{\infty}(\Omega;\Sym) @>\ccurl>>
C^{\infty}(\Omega;\Sym)
@>\div >> C^{\infty}(\Omega;\mathbb{R}^3) @>>> 0\\
 @.                   @VV{\Pi^W_h}V   @VV\Pi^{\Theta}_hV @VV\Pi^{\Sigma}_h V @VV{\Pi^V_h}V  @.\\
\mathbb{T}  \hookrightarrow @.  W_h @>\eps  >> \Theta_h  @>\ccurl >> \Sigma_h @>\div>>
 V_h @>>> 0
\end{CD}
\end{equation}
commutes.  The spaces $\Sigma_h$ and $V_h$, and the associated
interpolation operators $\Pi_h^\Sigma$ and $\Pi_h^V$, have been
introduced above, so it remains to define the spaces $W_h$ and
$\Theta_h$,
and the associated interpolation operators.

The discrete complex \eqref{cd-discrete} can be defined for all
polynomial levels, i.e., the two final spaces $\Sigma_h$ and $V_h$
can be taken as any of the pairs in the family introduced in
Section~\ref{higher-order}. However, in order to simplify the
discussion below we will only discuss the lowest order case introduced in
Section~\ref{lowest-order}, i.e.,
$\Sigma_h \subset H(\div,\Omega;\Sym)$ consists
of piecewise quartic matrix fields  with linear divergence,
while $V_h \in L^2(\Omega;\R^3)$ is composed of piecewise linears.

We will first describe the corresponding space $\Theta_h \subset
H(\ccurl,\Omega;\Sym)$. Locally on each tetrahedron this space
consists
of functions in $\mathcal{P}_6(K;\Sym)$, which is a space of dimension
504.
In order to specify the degrees of freedom on a tetrahedron $K \in
\mathcal T_h$ we
introduce the polynomial space
$$
E(f) = \{ \, v \in \mathcal{P}_7(f;\mathbb{R}) \, | \, v|_{\partial
  f}=0, \textstyle\int_e \partial_m v \ud s =0, \, e \in \Delta_1(f)
\, \}
$$
for each face $f \in \Delta_2(K)$.
The dimension of $E(f)$ is $12$. Also for each $e \in \Delta_1(K)$
define the operator $\Gamma_e : C^\infty(K;\Sym) \to
C^\infty(e;e^\perp)$ by
$\Gamma_e(S) =2 \partial_s Q_s (S s) - \grad_{e^{\perp}} s' S s$.
Note that if $f \in \Delta_2(K)$, with $n=n_f$, and $e \in \Delta_1(f)$, then
$n'\Gamma_e(S) = s' \Lambda_f(S)s$ on $e$.

The 504 degrees of freedom used to define the finite element space
$\Theta_h$ are the following:

\begin{enumerate}
\item $S$ and $ \ccurl S $ at each vertex,
48 degrees of freedom,
\item $\int_e s'S s \ v \ud s$, $v \in
  \mathcal{P}_4(e;\mathbb{R})$, $e \in \Delta_1(K)$,
30 degrees of freedom,
\item $\int_e \Gamma_e(S) \cdot v  \ud s$, $v \in
  \mathcal{P}_5(e;Q_s\mathbb{R}^3)$, $e \in \Delta_1(K)$, 72 degrees of freedom,
\item $\int_e Q_s (\ccurl S ): W \ud s$, $\, W \in
\mathcal{P}_2(e;Q_s\Sym)$, $e \in \Delta_1(K)$, 90 degrees of freedom,
\item $\int_e \rote P_s S Q_s \ud s$, $e \in \Delta_1(K) $, 6 degrees of freedom,
\item $\int_f Q_n S Q_n : \eps_f(v) \ud x_f$, $v \in
b_f\mathcal{P}_4(f;Q_n\mathbb{R}^3)$, $f \in \Delta_2(K),$
120 degrees of freedom,
\item $\int_e
(s'\Lambda_f(S) m)v
\ud s$, $v \in \P_1(e;\R) $ with $\int_e v \ud s =0$, $f \in \Delta_2(K)$, $e \in \Delta_1(f),$  12 degrees of freedom,
\item $\int_f \Lambda_f(S) \ud x_f$, $\,
f \in \Delta_2(K),$
12 degrees of freedom,
\item $\int_f \Lambda_f(S): \gradf \gradff v \ud x_f$, $\, v \in
  E(f)$, $f \in \Delta_2(K)$,
48 degrees of freedom,
\item $\int_K \ccurl S : T \ud x,
\, T \in \mathcal{M}_4(K)$, 6 degrees of freedom.
\item $\int_K S : \eps(v)\ud x,
\, v \in \mathcal{P}^0_7(K; \R^3)$, 60 degrees of freedom.
\end{enumerate}
These degrees of freedom are unisolvent for the space
$\mathcal{P}_6(K;\Sym)$.
This is in fact a consequence of the following result.

\begin{lemma}\label{theta-unisolvent}
Let $f \in \Delta_2(K)$ be fixed and assume that
$S \in \mathcal{P}_6(K;\Sym)$ with all the degrees of freedom
(1)--(9) associated to all $g \in \Delta(f)$
equal to zero. Then $Q_n S Q_n$ and $\Lambda_f(S)$ are identically zero
on the face $f$.
\end{lemma}

\begin{proof}
Let $f \in \Delta_2(K)$ be fixed, and asssume that
$S \in \mathcal{P}_6(K;\Sym)$ has all degrees of freedom
associated with the subsimplexes of $f$ equal to zero.
We start by observing that the degrees of freedom (1)--(4) implies
that
\begin{equation}\label{edge-prop}
s'S s \equiv 0, \quad Q_s \ccurl S \equiv 0, \quad
\Gamma_e(S) \equiv 0, \text{\quad on $e \in \Delta_1(f)$}.
\end{equation}
Next we will show that
\begin{equation}\label{rotrot=0}
\rotf \rotff Q_n S Q_n = P_n (\ccurl S) P_n =  0
\end{equation}
on $f$. Note that by \eqref{edge-prop} this quantity vanishes on
$\partial f$, i.e., $\rotf \rotff Q_n S Q_n \in b_f\P_1(f;\R P_n).$
Furthermore, using \eqref{rel09}, for each
$U \in \P_1(f;\R P_n)$
\[
\int_f \rotf \rotff Q_n S Q_n : U \ud x_f
= - \int_{\partial f} [\rotf(Ss)\cdot(Un) +
(Ss) \cdot\curlf(U n )] \ud s.
\]
Also, by \eqref{rot-div} we obtain
$\rotf (Ss) = -[\partial_s(s'C_nSs) + \partial_m (m'C_nSs)]n$.
However, the fact that $m'\Gamma_e(S) =0$ on $e \in \Delta_1(f)$
implies that $\partial_ms'Ss = 2 \partial_s m'Ss$ or
$\partial_mm'C_nSs = - 2 \partial_s s'C_nSs$.
Therefore, we can conclude that
$\rotf(S s) = \partial_s(s'C_nSs)$ on $e \in \Delta_1(f)$.
On the other hand, since $s'Ss =0$ on $e$ we have from \eqref{rel01}
that
\begin{align*}
(Ss) \cdot \curlf( Un) &= -(Ss) \cdot [C_n\grad (n'Un)]
\\
&= (C_nSs)\cdot \grad(n'Un) = (s'C_nSs)\partial_s(n'Un).
\end{align*}
Since $S$ is zero at each vertex of $f$ we can therefore conclude that
\begin{equation}\label{introtrot=0}
\int_f \rotf \rotff Q_n S Q_n : U \ud x_f
= - \int_{\partial f} \partial_s [(s'C_nSs)(n'Un)] \ud s =0
\end{equation}
for all $U \in \P_1(f;\R P_n)$, and hence
\eqref{rotrot=0} follows.
As a consequence, we obtain from \eqref{pol-elas-2d1} that $Q_nSQ_n =
\eps_f(u)$
for a suitable $u \in \P_7(f;Q_n\R^3)$, where $u$ can be chosen such
that $\int_e s'u \ud s =0$ for each edge
$e \in \Delta_1(f)$. Therefore, since $s'Ss = \partial_s(s'u) = 0$
on each $e \in \Delta_1(f)$, we conclude that $s'u = 0$ on $\partial
f$, and that $u=0$ at each vertex.
Furthermore, on each edge $e \in \Delta_1(f)$
\[
0= m'\Gamma_e(S) = 2\partial_s(m'Ss) - \partial_m(s'Ss)=
\partial_s^2(m'u),
\]
and therefore $u =0$ on $\partial f$.
Hence,
$u$ vanishes on $f$ by the degrees of freedom (6), and so does $Q_nSQ_n$.

Next we will show that $\rotf \Lambda_f(S) = 0$
on $f$. By \eqref{rot-Lambda} we have that
$\rotf \Lambda_f(S) = C_n(\ccurl S)P_n$, and
therefore \eqref{edge-prop} implies that $\rotf \Lambda_f(S)$ is zero
on $\partial f$. Hence, it is enough to show that 
\begin{equation}\label{introtLambda=0}
\int_f \rotf \Lambda_f(S): V \ud x_f = 0, \quad V \in \P_1(f; Q_n \Sym
P_n).
\end{equation}
Furthermore, by the degrees of freedom (8) we obtain that
\[
\int_f \rotf \Lambda_f(S): V \ud x_f
= \int_{\partial f} [\Lambda_f(S)s]\cdot (Vn)\ud s
= \int_{\partial f} [m'\Lambda_f(S)s]\cdot(m'Vn) \ud s,
\]
for any $V \in \P_1(f; Q_n \Sym P_n)$,
where the last identity holds since $s' \Lambda_f(S)s =
n'\Gamma_e(S)=0$ on $e \in \Delta_1(f)$.
From (7) above, this will be zero if we
can show that
$\int_e m'\Lambda_f(S)s \ud s$ vanishes for all $e \in \Delta_1(f)$.
However,
$m'\Lambda_f(S)s = \partial_s (m'Sn) + \partial_m (s'Sn)-\partial_n(s'Sm)$.
Therefore the degrees of freedom (5) imply that
$\int_e m'\Lambda_f(S)s \ud s = \int_e \partial_s (m'Sn) \ud s = 0$.
We therefore conclude that $\rotf \Lambda_f(S)$ is zero on $f$,
and from \eqref{pol-elas-2d2} we can conclude that $\Lambda_f(S)
= \gradf \gradff v$ for a suitable scalar field $v \in \P_7(f;\R)$,
where we assume that $v$ is chosen to be zero at each vertex.
Furthermore, the fact that
$0 = n'\Gamma_e(S) = s' \Lambda_f(S) s = \partial_s^2 v$
implies that $v \equiv 0$ on $\partial f$. In particular,
$\gradf v =0$ at each vertex.
Therefore, $v \in E(f)$  since (7) implies that
$0 = \int_e s(s'\Lambda_f(S)m) \ud s = \int_e s\partial_s
\partial_m v \ud s
= -\int_e  \partial_m v \ud s$,
and, as a consequence of (9), $v=0$, and so is $\Lambda_f(S)$.
\end{proof}

The lemma above implies that the set of functionals  (1)--(11)
are unisolvent for the space $\mathcal{P}_6(K;\Sym)$.
This follows since if all the functionals (1)--(11) are zero then
$S \in \mathcal N_6$ by this lemma, and hence $\ccurl S \in \mathcal M_4$. 
By the degrees of freedom (10) $\ccurl S = 0$, and by Lemma~\ref{problemma}
and the degrees of freedom (11), $S=0$.

The finite element space $\Theta_h$ is defined
as all functions which belongs to $\mathcal{P}_6(K;\Sym)$
for all $K \in \T_h$ and with the continuity conditions induced by the
degrees of freedom. Hence, as a consequence of
Lemma~\ref{theta-unisolvent}, the variables $Q_n S Q_n$ and
$\Lambda_f(S)$
are continuous for all $S \in \Theta_h$ and $f \in \Delta_2(\T_h)$,
and therefore Theorem~\ref{piecewise-smooth} implies that
$\Theta_h \subset H(\ccurl,\Omega;\Sym)$.
The degrees of freedom (1)--(11) also
defines a
canonical interpolation operator $\Pi_h^\Theta : C^\infty(\Omega;\Sym)
\to \Theta_h$ by requiring that it reproduces
all the functionals.

Finally, we need to describe the finite element space $W_h \subset
H^1(\Omega;\mathbb{R}^3)$. Locally on each $K \in \T_h$ this space is
taken to be
$\P_7(K;\R^3)$, which is a space of dimension 360. A vector field $w \in
\mathcal{P}_7(K;\mathbb{R}^3)$
is uniquely determined by:
\begin{enumerate}
\item the values of $w$ and its first order derivatives
at each vertex, 48 degrees of freedom,
\item $\int_e w \cdot q \ud s$, $q \in \mathcal{P}_3(e;\R^3)$,
$e \in \Delta_1(K)$, 72 degrees of freedom,
\item $\int_f \eps_f(Q_n w) : \eps_f(q) \ud x_f$, $q \in b_f\mathcal{P}_4(f;Q_n\R^3)$,
$f \in \Delta_2(K)$, 120 degrees of freedom,
\item $\int_e \partial_m(n'w) \ud s$, $f \in \Delta_2(K)$,
$e \in \Delta_1(f)$, 12 degrees of freedom
\item $\int_f \gradf \gradff (n'w) : \gradf \gradff (q) \ud x_f$,
$q \in E(f)$, $f\in\Delta_2(K)$, 48 degrees of freedom,
\item  $\int_K \eps (w) \cdot \eps (q) \ud x$ for all $q \in 
\mathcal{P}^0_7(K;\R^3)$,
60 degrees of freedom.
\end{enumerate}
It is rather straightforward to check that this set of functionals is
unisolvent for the space $\mathcal{P}_7(K;\mathbb{R}^3)$.
The finite element space $W_h$ is then defined as the set of functions
which are locally in $\mathcal{P}_7(K;\mathbb{R}^3)$ and with the
continuity induced by these degrees of freedom. Hence, all
$w \in W_h$ are continuous across faces $f$ in $\Delta_2(\T_h)$,
and they are $C^1$ at the vertices.
The associated interpolation operator $\Pi^W_h: C^\infty(\Omega;\R^3)
\to W_h$ will
reproduce all degrees of freedom given by (1)--(6).

By using the degrees of freedom for the discrete spaces, 
one may also verify that the  
diagram \eqref{cd-final} commutes. This verification is tedious, and
we will drop most of the details. However, we will illustrate what needs to
be done. For example, as part of the verification of the
relation $\Pi_h^\Sigma \circ \ccurl = \ccurl \circ \Pi_h^\Theta$ we
must show that for $T = (I - \Pi_h^\Theta)S, \, S \in
C^\infty(\Omega;\Sym)$,
we have (cf. degrees of freedom (3) for $\Sigma_h$)
\begin{align*}
\int_f P_n (\ccurl T)P_n : U \ud x_f &= \int_f \rotf \rotff Q_n T Q_n : U
\ud x_f =0,\\
\int_f C_n (\ccurl T) P_n : V \ud x_f &= \int_f \rotf \Lambda_f(T) : V
\ud x_f =0,
\end{align*}
for $U \in \P_1(f, \R P_n)$, $V \in \P_1(f, Q_n \Sym P_n)$, and all $f
\in \Delta_2(\T_h)$.
However, the verifications of these identities are almost identical 
to arguments leading to the formulas \eqref{introtrot=0} and 
\eqref{introtLambda=0} in the proof of 
Lemma~\ref{theta-unisolvent} above. The two identities above show that the
degrees of freedom (3) for $\Sigma_h$ are zero for functions of the form
$(\Pi_h^\Sigma \ccurl -  \ccurl \Pi_h^\Theta)S \in
\Sigma_h$. Similar, but simpler, arguments can be used to show that
the other degrees of freedom for $\Sigma_h$ are zero. Hence, the identity
$\Pi_h^\Sigma \circ \ccurl = \ccurl \circ \Pi_h^\Theta$ follows, and 
corresponding verifications can be done for the identity
$\Pi_h^\Theta \circ \eps = \eps\circ \Pi_h^W$.

\bibliographystyle{amsplain}
\bibliography{elas3dfamily}

\end{document}